\documentclass[11pt,twoside]{article}

\usepackage{float}

\newcommand{\HeadTitle}{Interpolation of Reciprocal Powers and Complete Exponential Bell Polynomials}

\newcommand{\HeadTitleTwo}{\begin{center}
\Large{\textit{Interpolation of Reciprocal Powers and Complete Exponential Bell Polynomials}}
\end{center}}

\usepackage[margin=0.8in]{geometry}
\usepackage{amsmath, amssymb, amsfonts, mathtools}
\usepackage{amsthm}
\usepackage{mathrsfs}
\usepackage{stmaryrd}
\usepackage{esint}
\usepackage{eqnalign}

\usepackage{graphicx}
\usepackage{enumitem}
\usepackage[most]{tcolorbox}
\usepackage{float}
\usepackage{listings}
\usepackage{array}
\usepackage{booktabs}
\usepackage{xcolor}
\usepackage{emptypage}
\usepackage{tensor}
\usepackage{tabularx}

\usepackage{fancyhdr}
\usepackage{hyperref}
\usepackage[nameinlink,noabbrev]{cleveref}

\usepackage{titlesec}
\usepackage[T1]{fontenc}
\usepackage[utf8]{inputenc}
\usepackage{lmodern}
\usepackage{titling}

\usepackage[backend=biber,style=numeric]{biblatex}
\providecommand{\HeadTitle}{} 
\providecommand{\HeadTitleTwo}{} 
\providecommand{\HeadAuthor}{Luc Ramsès TALLA WAFFO} 

\titleformat{\section}[block]
  {\normalfont\large\bfseries\itshape\centering}
  {§\thesection.}
  {1em}
  {}

\titleformat{\subsection}
  {\normalfont\normalsize\bfseries}
  {\thesubsection}
  {1em}
  {}

\newtheorem{theorem}{Theorem}[section]

\newtheorem{proposition}[theorem]{Proposition}
\newtheorem{corollary}[theorem]{Corollary}

\crefname{example}{example}{examples}
\Crefname{example}{Example}{Examples}

\crefname{corollary}{corollary}{corollaries}
\Crefname{corollary}{Corollary}{Corollaries}

\crefname{definition}{definition}{definitions}
\Crefname{definition}{Definition}{Definitions}

\crefname{remark}{remark}{remarks}
\Crefname{remark}{Remark}{Remarks}

\crefname{conjecture}{conjecture}{conjectures}
\Crefname{conjecture}{Conjecture}{Conjectures}

\crefname{lemma}{lemma}{lemmas}
\Crefname{lemma}{Lemma}{Lemmas}

\crefname{proposition}{proposition}{propositions}
\Crefname{proposition}{Proposition}{Propositions}

\crefname{theorem}{theorem}{theorems}
\Crefname{theorem}{Theorem}{Theorems}

\numberwithin{equation}{section}

\newcommand{\Q}{\mathbb{Q}}
\newcommand{\bell}[1]{\mathcal{B}_{#1}}

\usepackage{fontspec}
\begin{document}

\thispagestyle{fancy}

\vspace{0.2cm}

\begin{center}
\Large{\HeadTitleTwo}
\end{center}

\hspace{3cm}

\begin{center}
Luc Ramsès TALLA WAFFO \\
Technische Universität Darmstadt\\
Karolinenplatz 5, 64289 Darmstadt, Germany\\
ramses.talla@stud.tu-darmstadt.de\\
\vspace{0.5cm}
\today
\end{center}

\begin{abstract}
We study the polynomials $P_n^{(m)}$ of degree at most $n$ that interpolate the reciprocal power $x^{-m}$ at the consecutive positive integers $1,\ldots,n+1$. A formal inversion of the normalized nodal polynomial yields an explicit representation in terms of complete exponential Bell polynomials evaluated at generalized harmonic numbers. We derive recurrences in both the degree and the power index, an operator-theoretic representation, a Newton expansion, and generating functions involving the polylogarithm and the Gauss hypergeometric function. A weighted sum of the interpolation polynomials is expressed through the Hurwitz zeta function (or the digamma function when $m=1$), yielding absolutely convergent series for integer zeta values in terms of monomial coefficients. We extend the same interpolation family to Dirichlet beta values, obtaining a fixed-monomial-coefficient formula with explicit rational weights independent of the power index, as well as special-value and polylogarithmic representations. The same framework recovers Connon's Bell-polynomial identity and connects it with multiple zeta-star values. Finally, we prove that every zero lies in the open right half-plane, classify the real zeros and establish simplicity for $m=1$, and derive a local scaling limit governed by the normalized incomplete gamma function. The corresponding asymptotic expansions of the zeros involve ordinary zeta values. For fixed degree, we also establish eventual real-rootedness as the power index tends to infinity, with an explicit sufficient bound.
\end{abstract}

\paragraph{Notation.}
Throughout, $n\ge0$ and $m\ge1$ are integers, unless otherwise stated. The symbol $[x^k]F(x)$ denotes the coefficient of $x^k$ in a polynomial or formal power series $F(x)$. We write $\bell r$ for the complete exponential Bell polynomial of order $r$, with $\bell0=1$ and $0!=1$, and use the generalized harmonic numbers
\[
 H_N^{(s)}=\sum_{j=1}^{N}\frac1{j^s},\qquad H_N=H_N^{(1)}\qquad(N,s\ge1).
\]
The rising factorial is $(a)_k=a(a+1)\cdots(a+k-1)$, with $(a)_0=1$. Empty products are equal to $1$ and empty sums to $0$. We use $\zeta(s)$ and $\zeta(s,a)$ for the Riemann and Hurwitz zeta functions, respectively, $\Gamma$ for the gamma function, $\psi=\Gamma'/\Gamma$ for the digamma function, and $\gamma$ for the Euler--Mascheroni constant. In the zero analysis, $N=n+1$ and $L_n=H_N$. Additional notation is introduced where needed.

\section*{Introduction}\label{sec:introduction}
Polynomial interpolation at consecutive integers provides a natural setting in which algebraic identities, finite differences, and special functions can be studied together. For integers $n\ge0$ and $m\ge1$, let $P_n^{(m)}\in\Q[x]$ denote the unique polynomial of degree at most $n$ satisfying
\begin{equation}\label{eq:interp}
 P_n^{(m)}(j)=\frac1{j^m}\qquad(1\le j\le n+1).
\end{equation}

\vspace{0.2cm}

Although the Lagrange and Newton formulas give immediate representations of this polynomial, the reciprocal-power data admit a complementary description through the inverse of the normalized nodal polynomial. In particular, the coefficient of degree $r$ in the formal logarithm of this inverse is $H_{n+1}^{(r)}/r$. The exponential generating function of the complete Bell polynomials then gives an explicit formula that is uniform in the degree and the power index. The relevant formal and symmetric-function identities belong to the classical theory described in \cite{Comtet1974,Stanley1999,Macdonald1995}.

\vspace{0.2cm}

A second motivation comes from the appearance of finite differences in zeta-function theory. The series of Ser and Hasse, the work of Connon on Bell polynomials and nonlinear Euler sums, and Newton-series methods for zeta values provide natural points of contact with the present interpolation problem \cite{Ser1926,Hasse1930,Connon2008,FlajoletVepstas2008,Blagouchine2018}. Rather than specializing the interpolating polynomial at zero alone, we retain its variable and study a weighted transform of the entire family. This gives a coefficient-extraction identity for $\zeta(m+k+1)$, as well as a regularized gamma limit. The previously known zeta identities recovered in this way are identified explicitly in the text. We also retain the original polynomials to represent Dirichlet beta values, comparing the resulting coefficient formulas with classical Hurwitz and polylogarithmic transformations and existing series accelerations.

\vspace{0.2cm}

The third objective is to investigate the zeros of the interpolating polynomials. A Laplace-transform representation permits a half-plane localization argument, while the product formula provides a more precise classification for the first reciprocal power. In the large-degree regime, the harmonic scaling connects the local zero distribution with the normalized incomplete gamma function and gives asymptotic corrections involving zeta values. The classical zero theory for the special functions involved supplies the context for these results \cite{Sedletskii2004,Temme1995,DLMF}. We also consider the complementary limit of large power index at fixed degree.

\vspace{0.2cm}

The paper is organized as follows. \Cref{sec:algebraic} develops the algebraic reduction and formal inversion. \Cref{sec:bell} gives the explicit Bell-polynomial representations and special values. \Cref{sec:generating} establishes the recurrence relations, the operator representation, and the generating functions. \Cref{sec:zeta} studies the connections with zeta functions and Euler sums. \Cref{sec:dirichlet-beta} extends the interpolation formulas to the Dirichlet beta function. \Cref{sec:zeros} is devoted to zero localization, scaling limits, and eventual real-rootedness. The concluding remarks summarize the results and identify further questions.

\section{Algebraic reduction and formal inversion}\label{sec:algebraic}
\subsection{The nodal polynomial and the auxiliary quotient}
We shall prove, in particular, that the degree of $P_n^{(m)}$ is exactly $n$. Write
\[
 D_n(x)=\prod_{j=1}^{n+1}(x-j)
       =\sum_{r=0}^{n+1}c_{r,n}x^r,
 \qquad c_{0,n}=(-1)^{n+1}(n+1)!,
\]
with the convention $c_{r,n}=0$ for $r>n+1$. It is convenient to introduce the normalization
\[
 A_n=\frac{(-1)^n}{(n+1)!},\qquad
 G_n(x)=\prod_{j=1}^{n+1}\left(1-\frac{x}{j}\right)
       =-A_nD_n(x).
\]

\begin{proposition}[Algebraic reduction]\label{prop:reduction}
For every $n\ge0$ and $m\ge1$, there exists a unique polynomial
\[
 Q_n^{(m)}(x)=\sum_{r=0}^{m-1}b_{r,n,m}x^r
\]
such that
\begin{equation}\label{eq:factor}
 x^mP_n^{(m)}(x)-1=D_n(x)Q_n^{(m)}(x).
\end{equation}
This polynomial has degree exactly $m-1$, its constant coefficient is
\begin{equation}\label{eq:constant}
 b_{0,n,m}=\frac{(-1)^n}{(n+1)!},
\end{equation}
and all its coefficients are uniquely determined by the recurrence
\begin{equation}\label{eq:recursion}
 b_{r,n,m}=-\frac1{c_{0,n}}\sum_{k=0}^{r-1}
 b_{k,n,m}\,c_{r-k,n},\qquad 1\le r\le m-1.
\end{equation}
Moreover, $\deg P_n^{(m)}=n$.
\end{proposition}
\begin{proof}
By \eqref{eq:interp}, the polynomial $x^mP_n^{(m)}(x)-1$ vanishes at the $n+1$ distinct points $1,\ldots,n+1$; hence it is divisible by $D_n$. The quotient $Q_n^{(m)}$ is unique and initially has degree at most $m-1$. Comparing constant coefficients in \eqref{eq:factor} gives $-1=c_{0,n}b_{0,n,m}$, which proves \eqref{eq:constant}. For $1\le r<m$, the coefficient of $x^r$ on the left-hand side is zero. The Cauchy product therefore yields
\[
 0=\sum_{k=0}^{r}b_{k,n,m}c_{r-k,n},\qquad 1\le r<m,
\]
and, since $c_{0,n}\ne0$, this is precisely \eqref{eq:recursion}.

To establish the exact degrees, note that $D_n(x)=-A_n^{-1}G_n(x)$. Thus equation \eqref{eq:factor} becomes
\[
 G_n(x)\,\frac{Q_n^{(m)}(x)}{A_n}=1-x^mP_n^{(m)}(x).
\]
Consequently, in the formal power-series ring $\Q[[x]]$,
\begin{equation}\label{eq:inversecongruence}
 \frac{Q_n^{(m)}(x)}{A_n}\equiv G_n(x)^{-1}\pmod{x^m}.
\end{equation}
The formal expansion
\[
 G_n(x)^{-1}=\prod_{j=1}^{n+1}
 \left(\sum_{\ell=0}^{\infty}\frac{x^\ell}{j^\ell}\right)
\]
has a strictly positive coefficient in every degree. Hence the coefficient of $x^{m-1}$ in $Q_n^{(m)}$ is nonzero, so $\deg Q_n^{(m)}=m-1$. Comparing degrees in \eqref{eq:factor} now gives $\deg P_n^{(m)}=n$.
\end{proof}

\subsection{The first coefficients and Bell-polynomial notation}

The recurrence already shows that, for fixed $n$, the coefficient $b_{r,n,m}$ is independent of $m$ whenever $m>r$. To compute the first cases, set $h_s=H_{n+1}^{(s)}$ and write $b_{r,n,m}=A_n\beta_r$ for $0\le r<m$. Successive applications of \eqref{eq:recursion} give
\begin{align*}
 \beta_0={}&1,\\
 \beta_1={}&h_1,\\
 \beta_2={}&\frac{h_1^2+h_2}{2},\\
 \beta_3={}&\frac{h_1^3+3h_1h_2+2h_3}{6},\\
 \beta_4={}&\frac{h_1^4+6h_1^2h_2+3h_2^2+8h_1h_3+6h_4}{24},\\
 \beta_5={}&\frac{h_1^5+10h_1^3h_2+15h_1h_2^2+20h_1^2h_3
                   +20h_2h_3+30h_1h_4+24h_5}{120}.
\end{align*}
Thus the six polynomials corresponding to $m=1,\ldots,6$ are
\begin{align*}
 Q_n^{(1)}(x)&=A_n,\\
 Q_n^{(2)}(x)&=A_n(1+\beta_1x),\\
 Q_n^{(3)}(x)&=A_n(1+\beta_1x+\beta_2x^2),\\
 Q_n^{(4)}(x)&=A_n(1+\beta_1x+\beta_2x^2+\beta_3x^3),\\
 Q_n^{(5)}(x)&=A_n(1+\beta_1x+\beta_2x^2+\beta_3x^3+\beta_4x^4),\\
 Q_n^{(6)}(x)&=A_n(1+\beta_1x+\beta_2x^2+\beta_3x^3+\beta_4x^4
                    +\beta_5x^5).
\end{align*}
These coefficients have the structure of the \emph{complete exponential Bell polynomials}. We recall their defining exponential generating function \cite{Comtet1974,Leung2018}:
\begin{equation}\label{eq:bellgenerating}
 \exp\!\left(\sum_{s=1}^{\infty}y_s\frac{x^s}{s!}\right)
   =\sum_{r=0}^{\infty}\frac{\bell r(y_1,\ldots,y_r)}{r!}x^r.
\end{equation}
Equivalently, for $r\ge1$, the complete exponential Bell polynomial has the finite combinatorial expansion
\begin{equation}\label{eq:bellfinite}
 \bell r(y_1,\ldots,y_r)
 =r!\!\!\sum_{\substack{\ell_1,\ldots,\ell_r\ge0\\
                   \ell_1+2\ell_2+\cdots+r\ell_r=r}}
 \prod_{s=1}^{r}\frac1{\ell_s!}
 \left(\frac{y_s}{s!}\right)^{\ell_s}.
\end{equation}
The calculations above suggest the general identity
\[
 \beta_r=\frac1{r!}\bell r\big(h_1,1!h_2,2!h_3,\ldots,(r-1)!h_r\big)
\]
for every $r\ge0$, the argument list being empty when $r=0$. The next section proves this identity by formal power series, using an exponential structure standard in enumerative combinatorics \cite{Stanley1999}.

\section{Explicit Bell-polynomial representations and special values}\label{sec:bell}

\subsection{Closed formulas for the auxiliary and interpolation polynomials}

\begin{proposition}[Bell-polynomial representation of the auxiliary quotient]\label{prop:bell}
For all $n\ge0$ and $m\ge1$, the auxiliary polynomial is given explicitly by
\begin{equation}\label{eq:Qbell}
 \boxed{\displaystyle
 Q_n^{(m)}(x)=\frac{(-1)^n}{(n+1)!}
 \sum_{r=0}^{m-1}\frac{x^r}{r!}
 \bell r\big(H_{n+1}^{(1)},1!H_{n+1}^{(2)},\ldots,
               (r-1)!H_{n+1}^{(r)}\big).}
\end{equation}
\end{proposition}
\begin{proof}
Since $G_n(0)=1$, the polynomial $G_n$ is invertible in $\Q[[x]]$. For $F_n(x)=G_n(x)^{-1}$, formal logarithms give
\begin{align*}
 \log F_n(x)
 &=-\sum_{j=1}^{n+1}\log\left(1-\frac{x}{j}\right)\\
 &=\sum_{s=1}^{\infty}\frac{x^s}{s}
   \sum_{j=1}^{n+1}\frac1{j^s}
 =\sum_{s=1}^{\infty}H_{n+1}^{(s)}\frac{x^s}{s}.
\end{align*}
All logarithms and exponentials are interpreted formally; no analytic convergence assumption is required. Thus
\[
 F_n(x)=\exp\!\left(\sum_{s=1}^{\infty}
 (s-1)!H_{n+1}^{(s)}\frac{x^s}{s!}\right).
\]
Substituting $y_s=(s-1)!H_{n+1}^{(s)}$ into \eqref{eq:bellgenerating} yields
\[
 F_n(x)=\sum_{r=0}^{\infty}\frac{x^r}{r!}
 \bell r\big(H_{n+1}^{(1)},1!H_{n+1}^{(2)},\ldots,
              (r-1)!H_{n+1}^{(r)}\big).
\]
By \eqref{eq:inversecongruence}, $Q_n^{(m)}/A_n$ is exactly the truncation of this series through degree $m-1$. This proves \eqref{eq:Qbell}.
\end{proof}

Combining \eqref{eq:factor} and \eqref{eq:Qbell} gives, for $x\ne0$, the closed interpolation formula
\begin{equation}\label{eq:Pgeneral}
 P_n^{(m)}(x)=\frac1{x^m}\left[1-
 \prod_{j=1}^{n+1}\left(1-\frac{x}{j}\right)
 \sum_{r=0}^{m-1}\frac{x^r}{r!}
 \bell r\big(H_{n+1}^{(1)},1!H_{n+1}^{(2)},\ldots,
              (r-1)!H_{n+1}^{(r)}\big)\right].
\end{equation}
The apparent singularity at $x=0$ is removable because the right-hand side agrees with the polynomial $P_n^{(m)}$. Its value there is determined by the next coefficient in the generating function.

\subsection{The value at the origin}

\begin{proposition}[Constant coefficient]\label{prop:zero}
The value of the interpolating polynomial at zero is
\begin{equation}\label{eq:Pzero}
 \boxed{\displaystyle
 P_n^{(m)}(0)=\frac1{m!}\bell m\big(H_{n+1}^{(1)},1!H_{n+1}^{(2)},\ldots,
             (m-1)!H_{n+1}^{(m)}\big).}
\end{equation}
Equivalently,
\begin{equation}\label{eq:PzeroSum}
 P_n^{(m)}(0)=\sum_{\substack{\alpha_1,\ldots,\alpha_{n+1}\ge0\\
                     \alpha_1+\cdots+\alpha_{n+1}=m}}
 \frac1{1^{\alpha_1}2^{\alpha_2}\cdots(n+1)^{\alpha_{n+1}}}.
\end{equation}
\end{proposition}
\begin{proof}
Write $F_n(x)=G_n(x)^{-1}=\sum_{r\ge0}\beta_rx^r$. \Cref{prop:reduction,prop:bell} give
\[
 x^mP_n^{(m)}(x)=1-G_n(x)\sum_{r=0}^{m-1}\beta_rx^r
 =G_n(x)\sum_{r=m}^{\infty}\beta_rx^r.
\]
Since $G_n(0)=1$, the coefficient of $x^m$ on the right is $\beta_m$, whereas the same coefficient on the left is $P_n^{(m)}(0)$. Hence $P_n^{(m)}(0)=\beta_m$, and \eqref{eq:Pzero} follows from \eqref{eq:bellgenerating}. Expanding the finite product
\[
 F_n(x)=\prod_{j=1}^{n+1}\left(\sum_{\alpha_j=0}^{\infty}
           \frac{x^{\alpha_j}}{j^{\alpha_j}}\right)
\]
and comparing coefficients of $x^m$ proves \eqref{eq:PzeroSum}.
\end{proof}

For example, $P_n^{(1)}(0)=H_{n+1}$ and
$P_n^{(2)}(0)=\frac12\big(H_{n+1}^2+H_{n+1}^{(2)}\big)$.
Formula \eqref{eq:PzeroSum} also identifies this value as the complete homogeneous symmetric polynomial of degree $m$ in $1,1/2,\ldots,1/(n+1)$; the classical generating-function treatment of these symmetric polynomials can be found in \cite[Chapter I, \S2]{Macdonald1995}.

\subsection{Evaluation at the first exterior integer}

The first integer beyond the interpolation nodes, namely $n+2$, has an equally explicit description.

\begin{proposition}[Evaluation at the first exterior integer]\label{prop:next}
For all $n\ge0$ and $m\ge1$,
\begin{equation}\label{eq:PnextBell}
 \boxed{\displaystyle
 P_n^{(m)}(n+2)=\frac{1+(-1)^n\displaystyle\sum_{r=0}^{m-1}
       \frac{(n+2)^r}{r!}
       \bell r\big(H_{n+1}^{(1)},1!H_{n+1}^{(2)},\ldots,
               (r-1)!H_{n+1}^{(r)}\big)}{(n+2)^m}.}
\end{equation}
The same value has the alternating binomial representation
\begin{equation}\label{eq:PnextBinom}
 \boxed{\displaystyle
 P_n^{(m)}(n+2)=\sum_{j=1}^{n+1}
     (-1)^{n+1-j}\binom{n+1}{j-1}\frac1{j^m}.}
\end{equation}
In particular, $P_n^{(m)}(n+2)>0$ when $n$ is even. When $n$ is odd,
$P_n^{(1)}(n+2)=0$ and $P_n^{(m)}(n+2)<0$ for $m\ge2$.
\end{proposition}
\begin{proof}
At $x=n+2$,
\[
 D_n(n+2)=\prod_{j=1}^{n+1}(n+2-j)=(n+1)!.
\]
Substituting into \eqref{eq:factor} and applying \eqref{eq:Qbell} gives
\[
 (n+2)^mP_n^{(m)}(n+2)-1=(n+1)!\,Q_n^{(m)}(n+2),
\]
which proves \eqref{eq:PnextBell}. For the second expression, use Lagrange interpolation:
\[
 P_n^{(m)}(x)=\sum_{j=1}^{n+1}\frac1{j^m}
      \prod_{\substack{1\le i\le n+1\\i\ne j}}
        \frac{x-i}{j-i}.
\]
The Lagrange basis polynomials at $x=n+2$ satisfy
\begin{align*}
 \prod_{\substack{1\le i\le n+1\\i\ne j}}
        \frac{n+2-i}{j-i}
 &=(-1)^{n+1-j}\frac{(n+1)!}{(j-1)!(n+2-j)!}\\
 &=(-1)^{n+1-j}\binom{n+1}{j-1},
\end{align*}
which establishes \eqref{eq:PnextBinom}. Finally, the coefficients $\beta_r$ of $F_n$ are all strictly positive and $\beta_0=1$. The sum in \eqref{eq:PnextBell} is therefore equal to $1$ when $m=1$ and strictly greater than $1$ when $m\ge2$. The asserted sign properties follow.
\end{proof}

Thus the coefficients of the auxiliary polynomial $Q_n^{(m)}$ and the special values $P_n^{(m)}(0)$ and $P_n^{(m)}(n+2)$ are governed by the same complete exponential Bell polynomials. The underlying mechanism is the formal inversion of $\prod_{j=1}^{n+1}(1-x/j)$; it also explains the stability of the coefficients of $Q_n^{(m)}$ as $m$ increases. We next use the same formal inversion to derive recurrence relations and generating functions for the two indices of the interpolation family.

\section{Recurrences, spectral representation, and generating functions}\label{sec:generating}

\subsection{Recurrence in the power index}

We extend the family by setting $P_n^{(0)}(x)=1$ for all $n\ge0$.
For convenience, write
\begin{equation}\label{eq:beta-series}
 \frac{1}{G_n(t)}=\sum_{r=0}^{\infty}\beta_{n,r}t^r,
 \qquad \beta_{n,r}=P_n^{(r)}(0),
\end{equation}
where the second equality follows from \cref{prop:zero} for
$r\ge1$ and from the convention $\beta_{n,0}=1$ for $r=0$.
The series may be interpreted formally, or analytically for $|t|<1$.

\begin{proposition}[Recurrence with respect to $m$]
\label{prop:recurrence-m}
For every $n\ge0$ and $m\ge0$, the polynomials satisfy
\begin{equation}\label{eq:recurrence-m}
 \boxed{\displaystyle
 P_n^{(m+1)}(x)=\frac{P_n^{(m)}(x)-P_n^{(m)}(0)G_n(x)}{x}.}
\end{equation}
The quotient is a polynomial, and the initial condition is $P_n^{(0)}(x)=1$.
In particular,
\[
 P_n^{(1)}(x)=\frac{1-G_n(x)}{x}.
\]
\end{proposition}
\begin{proof}
By \eqref{eq:Pgeneral} and \eqref{eq:beta-series}, we have
\[
 x^mP_n^{(m)}(x)=1-G_n(x)\sum_{r=0}^{m-1}\beta_{n,r}x^r,
\]
which remains valid for $m=0$ if the sum is interpreted as empty.
Replacing $m$ by $m+1$ and subtracting the preceding identity yields
\[
 x^{m+1}P_n^{(m+1)}(x)-x^mP_n^{(m)}(x)
 =-G_n(x)\beta_{n,m}x^m.
\]
Cancelling the common factor $x^m$ gives
\begin{equation}\label{eq:auxiliary-recurrence}
 xP_n^{(m+1)}(x)=P_n^{(m)}(x)-P_n^{(m)}(0)G_n(x).
\end{equation}
Since $G_n(0)=1$, the numerator in \eqref{eq:recurrence-m} is
indeed divisible by $x$. The displayed initial case follows by taking $m=0$.
\end{proof}

\subsection{The interpolation operator and its spectrum}

\begin{corollary}[Operator-theoretic representation]
\label{cor:operator-representation}
Let $\mathbb{Q}_n[x]$ denote the vector space of
polynomials over $\mathbb{Q}$ of degree at most $n$.
Define the linear operator
\[
T_n:\mathbb{Q}_n[x]\longrightarrow\mathbb{Q}_n[x]
\]
by
\[
T_n(P)(x)
:=
\frac{P(x)-P(0)G_n(x)}{x}.
\]
The numerator vanishes at $x=0$ and has degree at most $n+1$,
so the quotient belongs to $\mathbb{Q}_n[x]$.

Then, for every integer $m\geq 1$,
\[
\boxed{
P_n^{(m)}=T_n^m(1).
}
\]

Furthermore, $T_n$ is diagonalizable over $\mathbb{Q}$,
with spectrum
\[
\sigma(T_n)
=
\left\{
1,\frac12,\ldots,\frac{1}{n+1}
\right\}.
\]
\end{corollary}

\begin{proof}
By \cref{prop:recurrence-m}, we have
\[
P_n^{(m+1)}=T_n(P_n^{(m)}).
\]
Moreover, the initial condition implies that
\[
P_n^{(1)}=T_n(1).
\]
The first assertion thus follows by induction on $m$.

Now, let $\{\ell_j\}_{j=1}^{n+1}$ denote the Lagrange
basis associated with the interpolation nodes
$1,2,\ldots,n+1$, namely
\[
\ell_j(x)
=
\prod_{\substack{1\leq k\leq n+1\\k\neq j}}
\frac{x-k}{j-k}.
\]

Since $G_n(k)=0$ for every $k\in\{1,\ldots,n+1\}$,
we obtain
\[
T_n(\ell_j)(k)
=
\frac{\ell_j(k)}{k}
=
\frac{\delta_{jk}}{k}
=
\frac{\delta_{jk}}{j}.
\]

Consequently, by the uniqueness of polynomial
interpolation,
\[
T_n(\ell_j)=\frac{1}{j}\ell_j.
\]

Thus, the Lagrange basis constitutes an eigenbasis
of $T_n$, and the asserted spectral characterization
follows.
\end{proof}

\subsection{Generating function in the power index}

\begin{proposition}[Generating function in the power index]
\label{prop:gf-m}
For every $n\ge0$, define
\[
 \mathcal{M}_n(x,t)=\sum_{m=1}^{\infty}P_n^{(m)}(x)t^{m-1}.
\]
Then, for $|t|<1$,
\begin{equation}\label{eq:gf-m}
 \boxed{\displaystyle
 \mathcal{M}_n(x,t)=\frac{G_n(t)-G_n(x)}{(x-t)G_n(t)}.}
\end{equation}
The apparent singularity at $x=t$ is removable. Equivalently,
\begin{equation}\label{eq:gf-m-zero}
 \sum_{m=0}^{\infty}P_n^{(m)}(x)t^m
 =\frac{xG_n(t)-tG_n(x)}{(x-t)G_n(t)}.
\end{equation}
\end{proposition}
\begin{proof}
The series defining $\mathcal{M}_n$ converges locally uniformly in $x$
for $|t|<1$: by the Lagrange representation, for each compact set
$K\subset\mathbb{C}$ the quantities $|P_n^{(m)}(x)|$ are bounded uniformly
for $x\in K$ and $m\ge1$.
Multiplying \eqref{eq:auxiliary-recurrence} by $t^m$ and summing over
$m\ge0$, we obtain, using \eqref{eq:beta-series},
\[
 x\mathcal{M}_n(x,t)=1+t\mathcal{M}_n(x,t)-\frac{G_n(x)}{G_n(t)}.
\]
Solving for $\mathcal{M}_n$ proves \eqref{eq:gf-m}. Since
$G_n(t)-G_n(x)$ is divisible by $x-t$, the apparent singularity is
removable. The identity also holds formally in $\Q[x][[t]]$,
with the divided difference interpreted as a polynomial.
Finally, \eqref{eq:gf-m-zero} follows from
$1+t\mathcal{M}_n(x,t)$.
\end{proof}

\begin{corollary}[Diagonal generating function]\label{cor:gf-diagonal}
For every $n\ge0$ and $|t|<1$,
\begin{equation}\label{eq:gf-diagonal}
 \sum_{m=1}^{\infty}P_n^{(m)}(t)t^{m-1}
 =\sum_{j=1}^{n+1}\frac{1}{j-t}
 =\sum_{r=0}^{\infty}H_{n+1}^{(r+1)}t^r.
\end{equation}
\end{corollary}
\begin{proof}
Taking the limit $x\to t$ in \eqref{eq:gf-m} gives
$\mathcal{M}_n(t,t)=-G_n'(t)/G_n(t)$.
Logarithmic differentiation of the product defining $G_n$ yields the
first equality in \eqref{eq:gf-diagonal}; expansion of each summand as a
geometric series gives the second.
\end{proof}

\subsection{Recurrence in the degree and Newton representation}

\begin{proposition}[Recurrence in the degree]
\label{prop:recurrence-n}
For all $n\ge1$ and $m\ge1$, the polynomials satisfy
\begin{equation}\label{eq:recurrence-n}
 \boxed{\displaystyle
 P_n^{(m)}(x)=P_{n-1}^{(m)}(x)
 +\frac{P_n^{(m-1)}(0)}{n+1}G_{n-1}(x),}
\end{equation}
with $P_0^{(m)}(x)=1$.
Equivalently, the correction term in \eqref{eq:recurrence-n} is
\[
 \frac{(-1)^n P_n^{(m-1)}(0)}{(n+1)!}\,D_{n-1}(x).
\]
\end{proposition}
\begin{proof}
Since
$G_n(u)=(1-u/(n+1))G_{n-1}(u)$, \cref{prop:gf-m} yields
\begin{align*}
 \mathcal{M}_n(x,t)-\mathcal{M}_{n-1}(x,t)
 &=\frac{1}{x-t}\left(
 \frac{G_{n-1}(x)}{G_{n-1}(t)}-\frac{G_n(x)}{G_n(t)}\right)\\
 &=\frac{G_{n-1}(x)}{(x-t)G_{n-1}(t)}
 \left(1-\frac{1-x/(n+1)}{1-t/(n+1)}\right)\\
 &=\frac{G_{n-1}(x)}{(n+1)G_n(t)}.
\end{align*}
The coefficient of $t^{m-1}$ on the right is
$G_{n-1}(x)\beta_{n,m-1}/(n+1)$ by \eqref{eq:beta-series}.
Comparing coefficients proves \eqref{eq:recurrence-n}. The initial condition
follows from the interpolation condition \eqref{eq:interp} when $n=0$.
\end{proof}

\begin{corollary}[Newton representation]\label{cor:newton}
For every $n\ge0$ and $m\ge1$,
\begin{equation}\label{eq:newton-representation}
 \boxed{\displaystyle
 P_n^{(m)}(x)=\sum_{k=0}^{n}
 \frac{\beta_{k,m-1}}{k+1}
 \prod_{j=1}^{k}\left(1-\frac{x}{j}\right).}
\end{equation}
Here the empty product is $1$. In particular, the leading coefficient is
\begin{equation}\label{eq:leading-coeff-newton}
 [x^n]P_n^{(m)}(x)=\frac{(-1)^n}{(n+1)!}\,\beta_{n,m-1}.
\end{equation}
\end{corollary}
\begin{proof}
Summing \eqref{eq:recurrence-n} over the degree index gives
\eqref{eq:newton-representation}, since the term with $k=0$ is $1$.
Only the term with $k=n$ contributes to the coefficient of $x^n$,
which proves \eqref{eq:leading-coeff-newton}.
\end{proof}

By \cref{prop:zero}, the Newton coefficients in
\eqref{eq:newton-representation} admit the explicit Bell-polynomial expression
\[
 \beta_{k,m-1}=\frac{1}{(m-1)!}\bell{m-1}\big(
 H_{k+1}^{(1)},1!H_{k+1}^{(2)},\ldots,(m-2)!H_{k+1}^{(m-1)}\big),
\]
where the right-hand side is interpreted as $1$ for $m=1$.
Unlike \eqref{eq:Pgeneral}, the Newton representation exhibits the
polynomial character of $P_n^{(m)}$ without a division by $x^m$.

\subsection{Generating function in the degree}

\begin{proposition}[Generating function in the degree]
\label{prop:gf-n}
For each integer $m\ge1$, let
\[
 \mathcal{N}_m(x,z)=\sum_{n=0}^{\infty}P_n^{(m)}(x)z^n.
\]
For $x\in\mathbb{C}$ and $|z|<1$, one has
\begin{equation}\label{eq:gf-n}
 \boxed{\displaystyle
 \mathcal{N}_m(x,z)=\frac{1}{(m-1)!(1-z)}
 \int_0^1(-\log u)^{m-1}(1-z+zu)^{x-1}\,du.}
\end{equation}
The complex power is defined using the principal logarithm.
\end{proposition}
\begin{proof}
For $k\ge0$ and $|s|<1$, the beta integral gives
\begin{align*}
 \int_0^1 u^{-s}(1-u)^k\,du
 &=\frac{k!}{(1-s)(2-s)\cdots(k+1-s)}\\
 &=\frac{1}{(k+1)G_k(s)}.
\end{align*}
Expanding $u^{-s}=\exp(-s\log u)$ in powers of $s$ and using
\eqref{eq:beta-series}, we deduce that
\begin{equation}\label{eq:beta-integral}
 \frac{\beta_{k,m-1}}{k+1}
 =\frac{1}{(m-1)!}\int_0^1(-\log u)^{m-1}(1-u)^k\,du.
\end{equation}
Termwise integration is justified on $|s|\le\rho<1$ by the integrable
majorant $u^{-\rho}$.
Since
\[
 \prod_{j=1}^{k}\left(1-\frac{x}{j}\right)=\frac{(1-x)_k}{k!},
\]
inserting \eqref{eq:beta-integral} into \eqref{eq:newton-representation}
gives
\[
 P_n^{(m)}(x)=\frac{1}{(m-1)!}\int_0^1(-\log u)^{m-1}
 \sum_{k=0}^{n}\frac{(1-x)_k}{k!}(1-u)^k\,du.
\]
For $|z|<1$, summing over $n$ and interchanging the sum and integral yields
\begin{align*}
 \mathcal{N}_m(x,z)
 &=\frac{1}{(m-1)!(1-z)}\int_0^1(-\log u)^{m-1}
 \sum_{k=0}^{\infty}\frac{(1-x)_k}{k!}[z(1-u)]^k\,du\\
 &=\frac{1}{(m-1)!(1-z)}\int_0^1(-\log u)^{m-1}
 (1-z+zu)^{x-1}\,du,
\end{align*}
where the last equality is the generalized binomial theorem.
For completeness, if $|x|\le R$, then
$|(1-x)_k/k!|\le\prod_{j=1}^{k}(1+R/j)$, which grows at most
polynomially in $k$. Consequently, the preceding interchanges are
absolutely and locally uniformly convergent for $|z|<1$;
$(-\log u)^{m-1}$ is integrable on $(0,1)$. Also,
$\operatorname{Re}(1-z+zu)>0$ on this domain, so the logarithm is
well defined.
\end{proof}

\subsection{Elementary and polylogarithmic specializations}

\begin{corollary}[Elementary and polylogarithmic specializations]
\label{cor:gf-n-special}
For $|z|<1$, the specialization at $m=1$ is
\begin{equation}\label{eq:gf-n-m1}
 \boxed{\displaystyle
 \mathcal{N}_1(x,z)=\frac{1-(1-z)^x}{xz(1-z)}.}
\end{equation}
Moreover, for every $m\ge1$,
\begin{equation}\label{eq:gf-n-polylog}
 \boxed{\displaystyle
 \mathcal{N}_m(0,z)=-\frac{1}{z(1-z)}\operatorname{Li}_m\left(\frac{z}{z-1}\right)},
\end{equation}
where $\operatorname{Li}_m$ denotes the polylogarithm on its principal
branch. All apparent singularities in these expressions at $x=0$ or
$z=0$ are interpreted by continuity.
\end{corollary}
\begin{proof}
Equation \eqref{eq:gf-n-m1} follows by direct integration in
\eqref{eq:gf-n} with $m=1$. For \eqref{eq:gf-n-polylog}, set $x=0$.
For $|z|<1/2$, expand
\[
 \frac{1}{1-z+zu}=\frac{1}{1-z}
 \sum_{k=0}^{\infty}\left(\frac{z}{z-1}\right)^k u^k.
\]
Using
\[
 \int_0^1u^k(-\log u)^{m-1}\,du=\frac{(m-1)!}{(k+1)^m},
\]
we obtain
\[
 \mathcal{N}_m(0,z)=\frac{1}{(1-z)^2}
 \sum_{k=0}^{\infty}\frac{(z/(z-1))^k}{(k+1)^m}.
\]
The defining series $\operatorname{Li}_m(w)=\sum_{j\ge1}w^j/j^m$
then gives \eqref{eq:gf-n-polylog}. The identity extends to $|z|<1$
by analyticity, since $w=z/(z-1)$ maps this disc into
$\operatorname{Re}(w)<1/2$, away from the principal branch cut of the
polylogarithm. At $z=0$ both generating functions have value $1$.
\end{proof}

In view of \cref{prop:zero}, equation \eqref{eq:gf-n-polylog}
can be restated as a generating-function identity for complete exponential
Bell polynomials evaluated at generalized harmonic numbers:
\begin{equation}\label{eq:bell-polylog}
\sum_{n=0}^{\infty}\frac{z^n}{m!}\bell m\big(
 H_{n+1}^{(1)},1!H_{n+1}^{(2)},\ldots,(m-1)!H_{n+1}^{(m)}\big)=-\frac{1}{z(1-z)}\operatorname{Li}_m\left(\frac{z}{z-1}\right).
\end{equation}
The identity \eqref{eq:gf-n-polylog} is also the index-shifted form of the multiple harmonic-sum generating function in \cite[Eqs.~(32)--(33)]{Boyadzhiev2009Euler}, by \eqref{eq:PzeroSum}. For instance, $m=1$ recovers the ordinary generating function
$-\log(1-z)/(z(1-z))$ of the shifted harmonic numbers $H_{n+1}$.

\subsection{The bivariate hypergeometric generating function}

\begin{proposition}[Bivariate generating function]\label{prop:gf-bivariate}
Define
\[
 \mathcal{F}(x,z,s)=\sum_{n=0}^{\infty}\sum_{m=1}^{\infty}
 P_n^{(m)}(x)z^ns^{m-1}.
\]
For $x\in\mathbb{C}$, $|z|<1$, and $|s|<1$, one has
\begin{align}
 \mathcal{F}(x,z,s)
 &=\frac{1}{1-z}\int_0^1 u^{-s}(1-z+zu)^{x-1}\,du
 \label{eq:gf-bivariate-integral}\\
 &=\boxed{\displaystyle\frac{{}_2F_1(1-x,1;2-s;z)}{(1-z)(1-s)}}.
 \label{eq:gf-bivariate}
\end{align}
Here ${}_2F_1$ is the Gauss hypergeometric function.
\end{proposition}
\begin{proof}
Summing \eqref{eq:gf-n} against $s^{m-1}$ gives
\eqref{eq:gf-bivariate-integral}, since
\[
 \sum_{m=1}^{\infty}\frac{(-\log u)^{m-1}}{(m-1)!}s^{m-1}=u^{-s}.
\]
The double series and the interchanges are absolutely convergent on
compact subsets of the stated domain. Indeed, the bound for
$|(1-x)_k/k!|$ used in \cref{prop:gf-n}, together with
\eqref{eq:beta-integral}, gives a majorant proportional to
$u^{-\rho}\sum_{k\ge0}|(1-x)_k/k!|r^k$ for $|z|\le r<1$ and
$|s|\le\rho<1$ (with the additional bounded factor $1/(1-r)$).
This is integrable in $u$.
After the substitution $v=1-u$, the integral becomes
\[
 \frac{1}{1-z}\int_0^1(1-v)^{-s}(1-zv)^{x-1}\,dv.
\]
Euler's integral representation
\[
 {}_2F_1(a,b;c;z)=\frac{\Gamma(c)}{\Gamma(b)\Gamma(c-b)}
 \int_0^1v^{b-1}(1-v)^{c-b-1}(1-zv)^{-a}\,dv
\]
applies with $a=1-x$, $b=1$, and $c=2-s$, since
$\operatorname{Re}(c-b)=1-\operatorname{Re}(s)>0$.
As $\Gamma(2-s)=(1-s)\Gamma(1-s)$, this proves
\eqref{eq:gf-bivariate}.
\end{proof}

The preceding results link the recurrence in the power index, the
recurrence in the degree index, and the Newton representation with a
single hypergeometric generating function. In particular,
\eqref{eq:gf-bivariate} simultaneously encodes the rational generating
functions in \eqref{eq:gf-m} and the integral generating functions in
\eqref{eq:gf-n}. These identities express the interpolation problem in
terms of established special-function and symmetric-function structures;
no claim of novelty for those classical structures is intended.

\section{Connections with zeta functions and Euler sums}\label{sec:zeta}

The preceding representations provide a natural link between reciprocal-power
interpolation and zeta functions. In particular, the Newton coefficients in
\eqref{eq:newton-representation} are the alternating finite differences of
reciprocal powers that occur in the classical series of Ser and Hasse
\cite{Ser1926,Hasse1930,Blagouchine2018}. Finite differences also arise in Newton-series
representations of the zeta function itself, as studied in
\cite{FlajoletVepstas2008}. Here we instead interpolate the reciprocal powers
and use the resulting polynomials as the summands of convergent zeta series.
We first recover a known Bell-polynomial identity and then retain the
interpolation variable to obtain identities for all monomial coefficients.

\subsection{Connon's identity and multiple zeta-star values}

For clarity, put
\[
 d_j^{(m)}:=\frac{\beta_{j,m-1}}{j+1}\qquad(j\ge0,\ m\ge1).
\]
Expanding $(1-u)^j$ in \eqref{eq:beta-integral} shows that
\begin{equation}\label{eq:zeta-difference}
 d_j^{(m)}=\sum_{\ell=0}^{j}(-1)^\ell\binom{j}{\ell}
 (\ell+1)^{-m}
 =\frac{1}{(m-1)!}\int_0^1(-\log u)^{m-1}(1-u)^j\,du.
\end{equation}
In particular, $0<d_j^{(m)}\le1$. Also, by \eqref{eq:PzeroSum},
\begin{equation}\label{eq:zeta-beta-bound}
 0<\beta_{j,r}\le H_{j+1}^{r}
 =O_r\big((1+\log(j+1))^r\big)\qquad(r\ge0).
\end{equation}
Indeed, the multinomial expansion of $H_{j+1}^{r}$ contains all the
terms in \eqref{eq:PzeroSum}, each with multiplicity at least one.

\begin{corollary}[Connon's Bell-polynomial identity]\label{cor:zeta-connon}
For every integer $m\ge2$,
\begin{equation}\label{eq:zeta-beta}
 \boxed{\displaystyle \zeta(m)=\frac{1}{m-1}
 \sum_{n=1}^{\infty}\frac{P_{n-1}^{(m-2)}(0)}{n^2}.}
\end{equation}
Equivalently, for $m\ge3$,
\begin{equation}\label{eq:zeta-connon}
 \boxed{\displaystyle
 \zeta(m)=\frac{1}{(m-1)!}\sum_{n=1}^{\infty}\frac{1}{n^2}
 \bell{m-2}\big(H_n^{(1)},1!H_n^{(2)},\ldots,
 (m-3)!H_n^{(m-2)}\big).}
\end{equation}
For $m=2$, the Bell factor in \eqref{eq:zeta-connon} is interpreted as
$\bell0=1$. Both series converge absolutely.
\end{corollary}
\begin{proof}
The classical Hasse representation \cite{Hasse1930,Blagouchine2018} gives
\begin{align*}
 (m-1)\zeta(m)
 &=\sum_{j=0}^{\infty}\frac{1}{j+1}
   \sum_{\ell=0}^{j}(-1)^\ell\binom{j}{\ell}(\ell+1)^{1-m}\\
 &=\sum_{j=0}^{\infty}\frac{d_j^{(m-1)}}{j+1}.
\end{align*}
For completeness, this instance of Hasse's formula also follows directly
from \eqref{eq:zeta-difference}: by Tonelli's theorem and the identity
$\sum_{j\ge0}(1-u)^j/(j+1)=-\log u/(1-u)$, the last series equals
\[
 \frac{1}{(m-2)!}\int_0^1\frac{(-\log u)^{m-1}}{1-u}\,du
 =(m-1)\zeta(m),
\]
using the standard integral for the Riemann zeta function
\cite[Eq.~25.5.1]{DLMF}. Now substitute
$d_j^{(m-1)}=\beta_{j,m-2}/(j+1)$, put $n=j+1$, and use
\eqref{eq:beta-series} to obtain \eqref{eq:zeta-beta}.
Finally, \eqref{eq:Pzero} gives \eqref{eq:zeta-connon} since
$(m-1)(m-2)!=(m-1)!$. Absolute convergence follows also from
\eqref{eq:zeta-beta-bound}.
\end{proof}

Identity \eqref{eq:zeta-connon} is not new. It is the specialization
$a=1$, $q=m-1$ of Connon's Hurwitz-zeta expansion
\cite[Abstract, p.~1]{Connon2008}, in which the factor $1/n^2$ is replaced
by $\Gamma(n)\Gamma(a)/(n\Gamma(n+a))$ and $H_n^{(s)}$ is replaced by
$\sum_{j=0}^{n-1}(a+j)^{-s}$. The derivation above recovers this known
specialization directly from our interpolation family. For $m=3,4,5$ it yields
\begin{align*}
 \zeta(3)&=\frac12\sum_{n=1}^{\infty}\frac{H_n}{n^2},\\
 \zeta(4)&=\frac16\sum_{n=1}^{\infty}\frac{H_n^2+H_n^{(2)}}{n^2},\\
 \zeta(5)&=\frac1{24}\sum_{n=1}^{\infty}
 \frac{H_n^3+3H_nH_n^{(2)}+2H_n^{(3)}}{n^2}.
\end{align*}
These are examples of nonlinear Euler sums; see also \cite{XuCai2018}.

The same identity has a multiple-zeta interpretation. For an admissible
index $(s_1,\ldots,s_r)$ of positive integers with $s_1\ge2$, define
\[
 \zeta^\star(s_1,\ldots,s_r)
 :=\sum_{n_1\ge\cdots\ge n_r\ge1}
 \frac{1}{n_1^{s_1}\cdots n_r^{s_r}}.
\]
By \eqref{eq:PzeroSum}, $\beta_{n-1,r}$ equals
$\sum_{n\ge k_1\ge\cdots\ge k_r\ge1}(k_1\cdots k_r)^{-1}$, with the
empty sum over tuples interpreted as $1$ when $r=0$. Thus
\eqref{eq:zeta-beta} recovers the height-one specialization of the
multiple-zeta-star sum formula:
\begin{equation}\label{eq:zeta-star}
 \zeta^\star(2,\underbrace{1,\ldots,1}_{m-2})=(m-1)\zeta(m).
\end{equation}
For the theory of multiple harmonic sums and the non-strict sum formula,
see \cite{Hoffman1992,OhnoOkuda2007}.

\subsection{A weighted Hurwitz-zeta transform}

We now retain $x$ rather than specializing it to zero.

\begin{proposition}[Hurwitz-zeta transform]\label{prop:zeta-transform}
For $m\ge1$ and $\operatorname{Re}x>-1$, the series
\[
 S_m(x):=\sum_{n=0}^{\infty}\frac{P_n^{(m)}(x)}{(n+1)(n+2)}
\]
converges absolutely and locally uniformly. On this half-plane,
\begin{equation}\label{eq:zeta-transform}
 \boxed{\displaystyle S_m(x)=
 \begin{cases}
 \dfrac{\zeta(m)-\zeta(m,1+x)}{x},&m\ge2,\\[5pt]
 \dfrac{\psi(1+x)+\gamma}{x},&m=1.
 \end{cases}}
\end{equation}
The expressions at $x=0$ are interpreted by analytic continuation, with
$S_m(0)=m\zeta(m+1)$ for all $m\ge1$.
\end{proposition}
\begin{proof}
Let $K$ be a compact subset of $\operatorname{Re}x>-1$, and put
$\sigma=\min_{x\in K}\operatorname{Re}x>-1$. Uniformly for $x\in K$,
\[
 \left|\frac{(1-x)_j}{j!}\right|=O_K((j+1)^{-\sigma}).
\]
To see this, separate finitely many factors in
$\prod_{\ell=1}^{j}(1-x/\ell)$, and use
$\log|1-x/\ell|=-\operatorname{Re}(x)/\ell+O_K(\ell^{-2})$
for the remaining factors. Since
$\sum_{n=j}^{\infty}((n+1)(n+2))^{-1}=(j+1)^{-1}$,
\eqref{eq:newton-representation} and \eqref{eq:zeta-beta-bound} imply
\begin{align*}
 \sum_{n\ge0}\frac{\sup_{x\in K}|P_n^{(m)}(x)|}{(n+1)(n+2)}
 &\le\sum_{j\ge0}\frac{d_j^{(m)}}{j+1}
       \sup_{x\in K}\left|\frac{(1-x)_j}{j!}\right|\\
 &\ll_{K,m}\sum_{j\ge0}\frac{(1+\log(j+1))^{m-1}}{(j+1)^{2+\sigma}}
 <\infty.
\end{align*}
This proves the asserted normal convergence and justifies the rearrangement
and, by \eqref{eq:zeta-difference}, the integral interchange below:
\begin{align*}
 S_m(x)&=\sum_{j=0}^{\infty}\frac{d_j^{(m)}}{j+1}
                 \frac{(1-x)_j}{j!}\\
 &=\frac1{(m-1)!}\int_0^1(-\log u)^{m-1}
       \sum_{j=0}^{\infty}\frac{(1-x)_j(1-u)^j}{(j+1)!}\,du\\
 &=\frac1{x(m-1)!}\int_0^1
         \frac{(-\log u)^{m-1}(1-u^x)}{1-u}\,du.
\end{align*}
Here we integrated the generalized binomial series to obtain
$\sum_{j\ge0}(1-x)_j v^j/(j+1)!=[1-(1-v)^x]/(xv)$ for $0<v<1$.
For $m\ge2$, the Hurwitz integral \cite[Eq.~25.11.25]{DLMF} identifies the
last integral with the difference in \eqref{eq:zeta-transform}. For $m=1$,
expanding $1/(1-u)$ gives
\[
 xS_1(x)=\sum_{\ell=1}^{\infty}\left(\frac1\ell-\frac1{\ell+x}\right)
 =\psi(1+x)+\gamma,
\]
the classical digamma expansion \cite[Eq.~5.7.6]{DLMF}.
Finally, taking the limit at zero in the integral gives
$S_m(0)=\int_0^1(-\log u)^m/(1-u)\,du/(m-1)!=m\zeta(m+1)$.
\end{proof}

\subsection{Zeta values from monomial coefficients}

\begin{corollary}[Zeta values from monomial coefficients]\label{cor:zeta-coefficients}
For all integers $m\ge1$ and $k\ge0$,
\begin{equation}\label{eq:zeta-coefficients}
 \boxed{\displaystyle
 \zeta(m+k+1)=\frac{(-1)^k}{\binom{m+k}{k+1}}
 \sum_{n=k}^{\infty}\frac{[x^k]P_n^{(m)}(x)}{(n+1)(n+2)}.}
\end{equation}
The series is absolutely convergent.
\end{corollary}
\begin{proof}
For $|x|<1$, expanding $1-u^x$ in the integral in the proof of
\cref{prop:zeta-transform} and integrating gives
\[
 S_m(x)=\sum_{k=0}^{\infty}(-1)^k\binom{m+k}{k+1}
              \zeta(m+k+1)x^k.
\]
Indeed the integral of $(-\log u)^{m+k}/(1-u)$ is
$(m+k)!\zeta(m+k+1)$. To justify termwise integration, fix $0<\rho<1$.
For $|x|\le\rho$ and $0<u<1$, with the quotient interpreted by continuity
at $x=0$, the absolute sum of the expansion of $(1-u^x)/x$ satisfies
\[
 \sum_{k=0}^{\infty}\frac{|x|^k(-\log u)^{k+1}}{(k+1)!}
 \le \frac{u^{-\rho}-1}{\rho}.
\]
Thus, apart from the constant factor $1/(m-1)!$, the absolute sum of
the expanded integrand is bounded by
\[
 \frac{(-\log u)^{m-1}(u^{-\rho}-1)}{\rho(1-u)}.
\]
This majorant is integrable: near $u=0$ it is
$O_{m,\rho}(u^{-\rho}(-\log u)^{m-1})$, while near $u=1$ it is
$O_{m,\rho}((1-u)^{m-1})$.
By the normal convergence in \cref{prop:zeta-transform} and Cauchy's
coefficient formula on any circle of radius $0<\rho<1$, we have
\[
 [x^k]S_m(x)=\sum_{n=k}^{\infty}\frac{[x^k]P_n^{(m)}(x)}{(n+1)(n+2)},
\]
and the series on the right converges absolutely. Comparing coefficients
proves the result.
\end{proof}

For example, the single family $P_n^{(2)}$ gives
\[
 \zeta(k+3)=\frac{(-1)^k}{k+2}\sum_{n=k}^{\infty}
       \frac{[x^k]P_n^{(2)}(x)}{(n+1)(n+2)}\qquad(k\ge0).
\]
The case $m=1$ is particularly simple: since $P_n^{(1)}=(1-G_n)/x$,
we have
\[
 (-1)^k [x^k]P_n^{(1)}(x)=e_{k+1}(1,1/2,\ldots,1/(n+1)),
\]
where $e_r$ is the elementary symmetric polynomial. More generally,
the Newton representation \eqref{eq:newton-representation} gives
\[
 (-1)^k[x^k]P_n^{(m)}(x)>0\qquad(0\le k\le n),
\]
since $d_j^{(m)}>0$ and the coefficient of $x^k$ in
$\prod_{\ell=1}^{j}(1-x/\ell)$ has sign $(-1)^k$ for $j\ge k$.
Thus \eqref{eq:zeta-coefficients} consists of nonnegative summands after
multiplication by $(-1)^k$, and it yields in particular
\begin{equation}\label{eq:zeta-elementary}
 \zeta(k+2)=\sum_{n=k}^{\infty}
 \frac{e_{k+1}(1,1/2,\ldots,1/(n+1))}{(n+1)(n+2)}.
\end{equation}
The underlying Hurwitz and digamma expansions are classical; the point
of \eqref{eq:zeta-coefficients} is to express them simultaneously in the
monomial coefficients of the interpolation family.

\subsection{Geometrically convergent zeta series}

The ordinary generating function \eqref{eq:gf-n-polylog} gives a second
connection with $\zeta(m)$ through the Euler transformation of the
alternating zeta series (see \cite{Sondow1994,DLMF}).

\begin{corollary}[Euler-transformed zeta series]\label{cor:zeta-accelerated}
For every integer $m\ge2$,
\begin{equation}\label{eq:zeta-accelerated}
 \boxed{\displaystyle
 \zeta(m)=\frac{1}{1-2^{1-m}}
 \sum_{n=0}^{\infty}\frac{P_n^{(m)}(0)}{2^{n+2}}.}
\end{equation}
Equivalently,
\begin{equation}\label{eq:zeta-accelerated-diff}
 \zeta(m)=\frac{1}{1-2^{1-m}}
 \sum_{n=0}^{\infty}\frac{P_n^{(m-1)}(0)}{2^{n+1}(n+1)}.
\end{equation}
The error after the first $N$ terms of \eqref{eq:zeta-accelerated-diff}
is nonnegative and at most $2^{-N}/(1-2^{1-m})$.
\end{corollary}
\begin{proof}
At $z=1/2$, \eqref{eq:gf-n-polylog} gives
$\mathcal{N}_m(0,1/2)=-4\operatorname{Li}_m(-1)
=4(1-2^{1-m})\zeta(m)$, proving \eqref{eq:zeta-accelerated}.
By \eqref{eq:newton-representation}, $P_n^{(m)}(0)=\sum_{j=0}^{n}d_j^{(m)}$.
Interchanging nonnegative sums and using
$\sum_{n=j}^{\infty}2^{-n-2}=2^{-j-1}$ gives
\eqref{eq:zeta-accelerated-diff}. Finally $0<d_j^{(m)}\le1$ by
\eqref{eq:zeta-difference}, so the stated error bound follows by summing
the geometric tail. These identities are the Euler-transformed eta series
in interpolation notation.
\end{proof}

\subsection{Spectral traces and the regularized gamma limit}

The operator $T_n$ from \cref{cor:operator-representation} gives the
finite spectral identities
\begin{equation}\label{eq:zeta-trace}
 \operatorname{tr}(T_n^m)=H_{n+1}^{(m)},\qquad
 \det(I-tT_n)=G_n(t).
\end{equation}
In particular, $\lim_{n\to\infty}\operatorname{tr}(T_n^m)=\zeta(m)$
for $m\ge2$. The determinant identity, together with
\eqref{eq:beta-series}, also yields a regularized generating function.

\begin{proposition}[Gamma limit]\label{prop:zeta-gamma}
Locally uniformly for $|t|<1$,
\begin{equation}\label{eq:zeta-gamma}
 \lim_{n\to\infty}e^{-tH_{n+1}}\sum_{r=0}^{\infty}
          P_n^{(r)}(0)t^r=e^{-\gamma t}\Gamma(1-t).
\end{equation}
Equivalently, the regularized logarithm converges to
$\sum_{r=2}^{\infty}\zeta(r)t^r/r$.
\end{proposition}
\begin{proof}
By the logarithmic expansion used in \cref{prop:bell},
\[
 \log G_n(t)^{-1}-tH_{n+1}
 =\sum_{r=2}^{\infty}H_{n+1}^{(r)}\frac{t^r}{r}.
\]
The bound $H_{n+1}^{(r)}\le\zeta(2)$ for $r\ge2$ justifies passage to the
limit locally uniformly on $|t|<1$. The classical expansion
$\log\Gamma(1-t)=\gamma t+\sum_{r\ge2}\zeta(r)t^r/r$
(see \cite[\S5.7(i)]{DLMF}) and exponentiation prove the claim.
\end{proof}

As a consequence, for $r\ge2$,
\begin{equation}\label{eq:zeta-regularized-bell}
 \lim_{n\to\infty}\sum_{j=0}^{r}\frac{(-H_{n+1})^j}{j!}
 P_n^{(r-j)}(0)
 =\frac{1}{r!}\bell r\big(0,1!\zeta(2),\ldots,(r-1)!\zeta(r)\big).
\end{equation}
Indeed, this follows by coefficient comparison in \eqref{eq:zeta-gamma}
and the defining Bell generating function \eqref{eq:bellgenerating}.
In contrast to the unregularized values $P_n^{(r)}(0)$, the combinations
on the left have finite limits. For instance, $r=2$ yields
$\lim_{n\to\infty}(2P_n^{(2)}(0)-(P_n^{(1)}(0))^2)=\zeta(2)$.
These limit identities are consequences of the classical gamma expansion,
not claims of new evaluations of zeta constants.

\section{Connections with the Dirichlet beta function}\label{sec:dirichlet-beta}

We retain the original interpolation nodes $1,\ldots,n+1$ and the
polynomials $P_n^{(m)}$ throughout. Let
\[
 \beta(s)=L(s,\chi_4)=\sum_{j=0}^{\infty}\frac{(-1)^j}{(2j+1)^s},
 \qquad \operatorname{Re}s>0,
\]
where $\chi_4$ is the nonprincipal character modulo $4$. The notation
$\beta(s)$ for the Dirichlet beta function is distinct from the coefficients
$\beta_{n,r}$ in \eqref{eq:beta-series}. We use the classical identities
\begin{align}
 \beta(m)&=4^{-m}\left[\zeta(m,1/4)-\zeta(m,3/4)\right],
      \qquad m\ge2,\label{eq:db-hurwitz}\\
 \beta(m)&=\frac{1}{(m-1)!}\int_0^1\frac{(-\log u)^{m-1}}{1+u^2}\,du
      =\operatorname{Im}\operatorname{Li}_m(i),\qquad m\ge1.
      \label{eq:db-integral}
\end{align}
The first identity is the character decomposition of the Hurwitz zeta
function (see \cite[\S25.15]{DLMF}). The integral follows by grouping the
alternating series, or by integrating the bounded partial geometric sums
of $\sum_{j\ge0}(-u^2)^j$. The last identity follows from the defining
polylogarithm series, with convergence at $i$ also for $m=1$.

\subsection{Special evaluations and positive coefficient series}

For convenience, throughout this section put
$\mathcal{A}_m(x):=xS_m(x)$, where $S_m$ is defined in
\cref{prop:zeta-transform} for $\operatorname{Re}x>-1$.

\begin{proposition}[Dirichlet beta values from the Hurwitz transform]
\label{prop:db-evaluations}
For every integer $m\ge1$,
\begin{equation}\label{eq:db-evaluations}
 \boxed{\displaystyle
 \beta(m)=\frac{1}{4^{m+1}}\sum_{n=0}^{\infty}
 \frac{3P_n^{(m)}(-3/4)-P_n^{(m)}(-1/4)}{(n+1)(n+2)}.}
\end{equation}
The series converges absolutely and all its summands are positive.
\end{proposition}
\begin{proof}
For $m\ge2$, \eqref{eq:zeta-transform}
and \eqref{eq:db-hurwitz} give
$4^m\beta(m)=\mathcal{A}_m(-1/4)-\mathcal{A}_m(-3/4)$.
For $m=1$ the same equality follows from the digamma expansion
$\psi(3/4)-\psi(1/4)=\sum_{j\ge0}[(j+1/4)^{-1}-(j+3/4)^{-1}]$.
Substituting the series for $S_m$ proves \eqref{eq:db-evaluations};
normal convergence follows from \cref{prop:zeta-transform}.
By the coefficient sign property following \eqref{eq:zeta-coefficients},
$P_n^{(m)}(-t)$ is positive and nondecreasing for $t\ge0$.
Consequently $3P_n^{(m)}(-3/4)-P_n^{(m)}(-1/4)>0$.
\end{proof}

\begin{corollary}[A positive shifted coefficient formula]\label{cor:db-shift}
For all $m\ge1$ and $k\ge0$, write $a_{n,j}^{(m)}=[x^j]P_n^{(m)}(x)$.
Then
\begin{align}
 \beta(m+k+1)
 &=\frac{(-1)^k}{4^{m+k+1}\binom{m+k}{k+1}}
   \sum_{n=k+1}^{\infty}\frac{[x^{k+1}]D_{n,m}^*(x)}{(n+1)(n+2)},
   \label{eq:db-shift}\\
 &=\frac{1}{4^{m+k+1}\binom{m+k}{k+1}}
   \sum_{n=k+1}^{\infty}\frac{1}{(n+1)(n+2)}
   \sum_{j=k+1}^{n}\binom{j+1}{k+1}
   \frac{3^{j-k}-1}{4^{j-k}}(-1)^ja_{n,j}^{(m)},
   \label{eq:db-shift-expanded}
\end{align}
where $D_{n,m}^*(x)=(x-3/4)P_n^{(m)}(x-3/4)
 -(x-1/4)P_n^{(m)}(x-1/4)$.
Both series converge absolutely, and the latter is a series of positive terms.
\end{corollary}
\begin{proof}
For $|x|<1/4$, define
\[
 F_m(x)=\sum_{j=0}^{\infty}\left[(j+1/4+x)^{-m}-(j+3/4+x)^{-m}\right].
\]
The series is normally convergent, including for $m=1$, by the difference
estimate $O_{m,K}((j+1)^{-m-1})$ on compact subsets of the disc.
By \eqref{eq:zeta-transform},
\[
 F_m(x)=\mathcal{A}_m(x-1/4)-\mathcal{A}_m(x-3/4)
 =-\sum_{n=0}^{\infty}\frac{D_{n,m}^*(x)}{(n+1)(n+2)}.
\]
The coefficient of $x^{k+1}$ in $F_m$ is
$(-1)^{k+1}\binom{m+k}{k+1}4^{m+k+1}\beta(m+k+1)$.
Cauchy's coefficient formula and normal convergence justify coefficient
comparison. The term of index $n=k$ vanishes because the highest-degree
terms of the two shifted polynomials cancel. The binomial theorem gives
\eqref{eq:db-shift-expanded}; positivity follows from
$(-1)^ja_{n,j}^{(m)}>0$ for $j\le n$.
\end{proof}

Expanding \eqref{eq:db-evaluations} at zero and using
\eqref{eq:zeta-coefficients}, with Tonelli's theorem for the nonnegative
terms, also yields
\begin{equation}\label{eq:db-zeta-series}
 \beta(m)=\sum_{k=0}^{\infty}\binom{m+k}{k+1}
       \frac{3^{k+1}-1}{4^{m+k+1}}\zeta(m+k+1),\qquad m\ge1.
\end{equation}
This is equivalently the classical Taylor expansion of the Hurwitz
zeta difference about the second argument $1$ (or of the digamma
difference for $m=1$), rather than an independent evaluation of beta values.
For a related and faster approach using reflection and horizontal shifts,
see \cite[\S7.1--7.2]{Languasco2023}.

\subsection{Beta values from a fixed monomial coefficient}

We next retain exactly the coefficient $[x^k]P_n^{(m)}$ in
\eqref{eq:zeta-coefficients}, changing only its weight. Put
\[
 e_{n,k}=e_k(1,1/2,\ldots,1/n)=[t^k]\prod_{j=1}^{n}(1+t/j),
\]
with $e_{n,0}=1$ and $e_{n,k}=0$ for $k>n$. For $n\ge k$ we have
$e_{n,k}>0$. Define rational numbers $q_r$ by
\begin{equation}\label{eq:db-q}
 q(v):=\sum_{r\ge0}q_rv^r=\frac{v}{2-2v+v^2}.
\end{equation}
In particular $q_0=0$, $q_1=1/2$, and
$q_r=q_{r-1}-q_{r-2}/2$ for $r\ge2$. Define for $n\ge k$
\begin{equation}\label{eq:db-weights}
 W_{n,k}=\frac{1}{e_{n,k}}\sum_{j=k}^{n}\frac{q_{n-j}e_{j,k}}{j+1},
 \qquad \omega_{n,k}=W_{n,k}-W_{n+1,k}.
\end{equation}
These rational weights depend on $n,k$, but not on $m$.

\begin{proposition}[Dirichlet beta values from monomial coefficients]
\label{prop:db-coefficients}
For all integers $m\ge1$ and $k\ge0$,
\begin{equation}\label{eq:db-coefficients}
 \boxed{\displaystyle
 \beta(m+k+1)=\frac{(-1)^k}{\binom{m+k}{k+1}}
     \sum_{n=k}^{\infty}\omega_{n,k}[x^k]P_n^{(m)}(x).}
\end{equation}
The series converges absolutely. For fixed $k$,
$W_{n,k}=(n+1)^{-1}+O_k((n+1)^{-2})$ and
$\omega_{n,k}=O_k((n+1)^{-2})$. For fixed $m,k$, the error after
truncating \eqref{eq:db-coefficients} at $n=N$ is
\[
 O_{m,k}\left(\frac{(1+\log(N+2))^{m+k}}{N+1}\right)
 \qquad(N\to\infty).
\]
\end{proposition}
\begin{proof}
Write $c_{n,k}^{(m)}=(-1)^k[x^k]P_n^{(m)}(x)$ and set
$c_{-1,k}^{(m)}=0$. The Newton representation
\eqref{eq:newton-representation} gives
\begin{equation}\label{eq:db-coeff-diff}
 c_{n,k}^{(m)}-c_{n-1,k}^{(m)}=e_{n,k}d_n^{(m)}.
\end{equation}
The generalized binomial theorem gives
$\sum_{n\ge0}\prod_{j=1}^n(1+t/j)v^n=(1-v)^{-1-t}$.
Integrating with respect to $v$, dividing by $v$, and extracting the
coefficient of $t^k$ shows that, for $|v|<1$,
\begin{equation}\label{eq:db-e-gf}
 \sum_{n=k}^{\infty}\frac{e_{n,k}}{n+1}v^n
 =\frac{(-\log(1-v))^{k+1}}{(k+1)!\,v}.
\end{equation}
The right-hand side is interpreted by continuity at $v=0$. Multiplying
by \eqref{eq:db-q} and comparing coefficients gives
\begin{equation}\label{eq:db-weight-gf}
 \sum_{n=k}^{\infty} W_{n,k}e_{n,k}v^n
 =\frac{(-\log(1-v))^{k+1}}{(k+1)!\,[1+(1-v)^2]}.
\end{equation}

We justify the estimates needed for integration and summation by parts.
The poles of $q(v)$ are $1\pm i$, so $q_r=O(2^{-r/2})$ and
$\sum_{r\ge0}q_r=q(1)=1$. For fixed $k\ge1$ and $n\ge k$, the elementary
symmetric functions satisfy $e_{n,k-1}\le k^2 e_{n,k}$: in the identity
$ke_{n,k}=\sum_{|I|=k-1}(\prod_{j\in I}j^{-1})\sum_{j\notin I}j^{-1}$,
with $I\subseteq\{1,\ldots,n\}$, at least one of the first $k$ indices
is absent from each $I$. Using
$e_{n,k}-e_{n-1,k}=e_{n-1,k-1}/n$, we obtain uniformly for
$0\le r\le n/2$ and large $n$
\[
 0\le 1-\frac{e_{n-r,k}}{e_{n,k}}
 \le k^2\sum_{j=n-r+1}^{n}\frac1j=O_k(r/n).
\]
For $k=0$ the ratio is identically one. Hence for $r\le n/2$,
\[
 \frac{e_{n-r,k}}{(n-r+1)e_{n,k}}
 =\frac1{n+1}+O_k\left(\frac{r}{(n+1)^2}\right).
\]
In the convolution defining $W_{n,k}$, the part $r>n/2$ is exponentially
small, since $0<e_{n-r,k}/e_{n,k}\le1$ for $r\le n-k$.
The finite first absolute moment of $(q_r)$ now gives
$W_{n,k}=(n+1)^{-1}+O_k((n+1)^{-2})$ and the asserted bound on $\omega_{n,k}$.

Furthermore $e_{n,k}\le H_n^k/k!$ for $n\ge1$, while
\eqref{eq:zeta-beta-bound} gives
$d_n^{(m)}=O_m((1+\log(n+1))^{m-1}/(n+1))$.
Thus $\sum_{n\ge k}|W_{n,k}|e_{n,k}d_n^{(m)}<\infty$, and
\eqref{eq:db-coeff-diff} implies
$c_{n,k}^{(m)}=O_{m,k}((1+\log(n+1))^{m+k})$.
With $v=1-u$ in \eqref{eq:db-weight-gf}, we may integrate termwise
against $(-\log u)^{m-1}/(m-1)!$ using
\eqref{eq:zeta-difference} and the preceding absolute bound. We obtain
\begin{align*}
 \sum_{n=k}^{\infty}W_{n,k}e_{n,k}d_n^{(m)}
 &=\frac1{(m-1)!(k+1)!}\int_0^1\frac{(-\log u)^{m+k}}{1+u^2}\,du\\
 &=\binom{m+k}{k+1}\beta(m+k+1).
\end{align*}
Finally, finite summation by parts in \eqref{eq:db-coeff-diff} gives
\[
 \sum_{n=k}^{N}W_{n,k}e_{n,k}d_n^{(m)}
 =\sum_{n=k}^{N}\omega_{n,k}c_{n,k}^{(m)}
    +W_{N+1,k}c_{N,k}^{(m)}.
\]
The last term tends to zero, and the series on the right converges
absolutely by the bounds above. This proves \eqref{eq:db-coefficients}.
Moreover, its absolute tail is bounded, up to a constant depending only
on $m,k$, by
\[
 \sum_{n>N}\frac{(1+\log(n+1))^{m+k}}{(n+1)^2}
 =O_{m,k}\left(\frac{(1+\log(N+2))^{m+k}}{N+1}\right),
\]
where the estimate follows from the integral test (for sufficiently
large $N$) and integration by parts. The constant binomial normalization
in \eqref{eq:db-coefficients} does not change this order.
\end{proof}

For comparison, the same proof with the multiplier $q(v)=1$ gives
$W_{n,k}=1/(n+1)$, $\omega_{n,k}=1/((n+1)(n+2))$, and recovers
\eqref{eq:zeta-coefficients}: the integral kernel is then $1/(1-u)$
in place of $1/(1+u^2)$. The weights in \eqref{eq:db-coefficients}
need not be positive. For instance, $\omega_{0,0}=-1/2$.
When $m=1$ and $k=0$, we obtain the particularly simple example
\[
 \beta(2)=\sum_{n=0}^{\infty}\omega_{n,0}H_{n+1}.
\]

\subsection{Geometrically convergent series from constant coefficients}

\begin{corollary}[Polylogarithmic beta series]\label{cor:db-polylog}
For every integer $m\ge1$,
\begin{equation}\label{eq:db-polylog}
 \boxed{\displaystyle
 \beta(m)=\frac12\sum_{n=0}^{\infty}
 2^{-n/2}\sin(n\pi/4)\,P_n^{(m)}(0).}
\end{equation}
Equivalently,
\begin{equation}\label{eq:db-blocks}
 \beta(m)=\sum_{j=0}^{\infty}\frac{(-1)^j}{4^{j+1}}
 \left[P_{4j+1}^{(m)}(0)+P_{4j+2}^{(m)}(0)
       +\tfrac12P_{4j+3}^{(m)}(0)\right].
\end{equation}
Both are absolutely convergent series with rational terms. The error
of truncating \eqref{eq:db-polylog} after $n=N$ is
$O_m(2^{-N/2}(1+\log(N+2))^m)$. In addition, if
$a_n=\tfrac12 2^{-n/2}\sin(n\pi/4)P_n^{(m)}(0)$ denotes its summand, then
\[
 \limsup_{n\to\infty}|a_n|^{1/n}=2^{-1/2}.
\]
\end{corollary}
\begin{proof}
Take $z=(1-i)/2$ in \eqref{eq:gf-n-polylog}. Then
$z/(z-1)=i$, $z(1-z)=1/2$ and $|z|=2^{-1/2}<1$.
Taking imaginary parts of
$\operatorname{Li}_m(i)=-\tfrac12\sum_{n\ge0}P_n^{(m)}(0)z^n$
proves \eqref{eq:db-polylog}. Grouping the weights in blocks of four
proves \eqref{eq:db-blocks}. The bound on $P_n^{(m)}(0)=\beta_{n,m}$
in \eqref{eq:zeta-beta-bound} yields the stated tail estimate, by summing
a geometric sequence times a fixed power of a logarithm. The same bound
implies $\limsup_{n\to\infty}|a_n|^{1/n}\le2^{-1/2}$. Conversely,
\eqref{eq:PzeroSum} gives $P_n^{(m)}(0)\ge1$ (the term with all $m$
powers assigned to the first node). Along $n\equiv2\pmod4$ one has
$|\sin(n\pi/4)|=1$, giving the reverse inequality.
\end{proof}

Substitution of \eqref{eq:Pzero} yields the Bell-polynomial form
\begin{equation}\label{eq:db-bell}
 \beta(m)=\frac{1}{2m!}\sum_{n=1}^{\infty}2^{-n/2}\sin(n\pi/4)
 \bell m\left(H_{n+1}^{(1)},1!H_{n+1}^{(2)},\ldots,
 (m-1)!H_{n+1}^{(m)}\right).
\end{equation}
In particular, for Catalan's constant $G=\beta(2)$,
\[
 G=\frac14\sum_{n=1}^{\infty}2^{-n/2}\sin(n\pi/4)
       \left(H_{n+1}^2+H_{n+1}^{(2)}\right).
\]

\subsection{Relation to earlier work}
The special-function identities used here are classical. More specifically,
Boyadzhiev's equations (32)--(33) in \cite{Boyadzhiev2009Euler}
are exactly \eqref{eq:gf-n-polylog} after shifting the index and
identifying his weakly ordered reciprocal sums with \eqref{eq:PzeroSum}.
Thus \eqref{eq:db-polylog} and \eqref{eq:db-bell} are direct
specializations of an established generating function, not independent
claims of novelty. The Hurwitz evaluations and the Taylor expansion
\eqref{eq:db-zeta-series} likewise use standard special-function relations.
Connon \cite{Connon2010Beta} obtained rapidly converging series for
Catalan's and Ap\'ery's constants. Languasco \cite[\S7.1--7.2]{Languasco2023}
uses expansions of Hurwitz zeta values to obtain computational formulas
for beta values. Hessami Pilehrood and Hessami Pilehrood
\cite{HessamiPilehrood2010Beta} obtained Ap\'ery-like series, including
one for Catalan's constant with geometric ratio $2^{-10}$, considerably
smaller than the exponential decay factor $2^{-1/2}$ per ungrouped
term in \eqref{eq:db-polylog}. The latter factor is understood in the
root-test sense established in \cref{cor:db-polylog}, not as a ratio of
successive summands, which is undefined whenever the denominator
vanishes. In contrast, the fixed-coefficient series
\eqref{eq:db-coefficients} has an algebraic-logarithmic upper bound on
its truncation error and is not intended as a series acceleration.
The recent preprint \cite[Prop.~3.3]{Bhandari2026Beta} also uses complete
Bell polynomials in Euler-transformed beta series, with harmonic sums
on shifted odd nodes, and studies reciprocal powers of beta values.

The purpose of the present extension is not to outperform those
accelerations or to claim new evaluations of classical constants. Its
specific feature is the simultaneous treatment of zeta and beta values
using the same consecutive-integer interpolation family. In particular,
\eqref{eq:db-coefficients} retains a single monomial coefficient of
$P_n^{(m)}$ and the binomial normalization of \eqref{eq:zeta-coefficients},
with explicit rational weights independent of $m$. This precise weighted
interpolation formulation is not identified in the references just discussed;
we make no exhaustive priority claim.

\section{Zero localization and asymptotic behaviour}\label{sec:zeros}

We now study the zeros of the interpolation polynomials themselves. The
product representation \eqref{eq:Pgeneral} leads to a half-plane localization
valid for every positive integer $m$, whereas the case $m=1$ admits a
complete classification of the real zeros. We then connect the local
asymptotic behaviour of the zeros to the incomplete gamma function and to
ordinary zeta values. Throughout this section, $n\ge1$, unless otherwise stated; we use the abbreviations $N=n+1$ and $L_n=H_N$ introduced above.

\subsection{Half-plane localization via Laplace transforms}

The following argument uses a Laplace transform with a positive decreasing
kernel. For related classical results on zero localization for Laplace
transforms, see \cite{Sedletskii2004}. We give a direct proof adapted to
our polynomial family.

\begin{proposition}[Half-plane localization]\label{prop:zeros-half-plane}
For all integers $n\ge1$ and $m\ge1$, every zero $\rho$ of $P_n^{(m)}$
satisfies
\begin{equation}\label{eq:zeros-half-plane}
 \boxed{\operatorname{Re}\rho>0.}
\end{equation}
Equivalently, the polynomial $P_n^{(m)}(-z)$ is Hurwitz stable.
\end{proposition}
\begin{proof}
For $t>0$ put $f_N(t)=Ne^{-t}(1-e^{-t})^{N-1}$. This is a strictly
positive probability density. Substitution $u=e^{-t}$ and the beta
integral give, for $\operatorname{Re}z>-1$,
\begin{equation}\label{eq:zeros-density}
 \int_0^\infty e^{-zt}f_N(t)\,dt
 =\frac{\Gamma(N+1)\Gamma(1+z)}{\Gamma(N+1+z)}
 =\frac1{G_n(-z)}.
\end{equation}
In particular, by \eqref{eq:beta-series},
$\displaystyle\int_0^\infty t^r f_N(t)\,dt=r!\beta_{n,r}$ for $r\ge0$.
Define
\[
 K_{n,m}(t)=\frac1{(m-1)!}
 \int_t^\infty (s-t)^{m-1}f_N(s)\,ds\qquad(t\ge0).
\]
Taylor's formula with integral remainder, applied to $e^{-zs}$, gives
\[
 e^{-zs}-\sum_{r=0}^{m-1}\frac{(-zs)^r}{r!}
 =\frac{(-z)^m}{(m-1)!}\int_0^s(s-t)^{m-1}e^{-zt}\,dt.
\]
Integrating this equality against $f_N(s)$ and using
\eqref{eq:Pgeneral}, \eqref{eq:beta-series} and \eqref{eq:zeros-density}
we obtain
\begin{equation}\label{eq:zeros-laplace}
 \frac{P_n^{(m)}(-z)}{G_n(-z)}
 =\int_0^\infty e^{-zt}K_{n,m}(t)\,dt
 \qquad(\operatorname{Re}z\ge0).
\end{equation}
Fubini's theorem applies since the absolute value of the integrated
remainder, after division by $(-z)^m$, is bounded by $s^m/m!$ on this
half-plane. The identity at $z=0$ follows by continuity.

The kernel is positive and tends to zero at infinity. Moreover,
$-K_{n,1}'=f_N$, and $-K_{n,m}'=K_{n,m-1}$ for $m\ge2$; in particular
$-K_{n,m}'(t)>0$ for $t>0$. Integration by parts gives
\[
 z\frac{P_n^{(m)}(-z)}{G_n(-z)}
 =K_{n,m}(0)-\int_0^\infty e^{-zt}(-K_{n,m}'(t))\,dt.
\]
If $\operatorname{Re}z>0$, the absolute value of the last integral is
strictly less than $K_{n,m}(0)$. If $\operatorname{Re}z=0$ and $z\ne0$,
the same strict inequality follows from the equality condition in the
triangle inequality: the phase of $e^{-zt}$ is not constant on the
support of the strictly positive density $-K_{n,m}'$. In either case
$P_n^{(m)}(-z)\ne0$, since $G_n(-z)\ne0$. Finally
$P_n^{(m)}(0)=\beta_{n,m}>0$ by \eqref{eq:PzeroSum}.
\end{proof}

\subsection{Simplicity and real zeros for the first reciprocal power}

\begin{proposition}[Simple zeros and a vertical strip]\label{prop:zeros-m1}
For $m=1$, all $n$ zeros of $P_n^{(1)}$ are simple and satisfy
\begin{equation}\label{eq:zeros-strip}
 0<\operatorname{Re}\rho\le n+2.
\end{equation}
If $n$ is even, there are no real zeros. If $n$ is odd, the unique
real zero is $n+2$, and it is the only zero on the line
$\operatorname{Re}z=n+2$.
\end{proposition}
\begin{proof}
By \eqref{eq:Pgeneral}, $P_n^{(1)}(z)=(1-G_n(z))/z$, with a removable
singularity at zero. Thus its zeros are precisely the nonzero solutions
of $G_n(z)=1$. The product satisfies
\begin{equation}\label{eq:zeros-reflection}
 G_n(N+1-z)=(-1)^NG_n(z).
\end{equation}
For $0<x<1$, clearly $0<G_n(x)<1$. For $k<x<k+1$ with
$1\le k<N$, estimating separately the factors with $j\le k$ and $j>k$
yields
\[
 |G_n(x)|<\frac{k!(N-k)!}{N!}=\binom Nk^{-1}<1.
\]
At the integer nodes the product vanishes, and \eqref{eq:zeros-reflection}
handles $N<x<N+1$. Consequently $|G_n(x)|<1$ on $(0,N+1)$.
At $x=N+1$, $G_n(N+1)=(-1)^N$, and for $x>N+1$ its modulus is
strictly greater than one. Negative real zeros are excluded by
\cref{prop:zeros-half-plane}. This proves the real-zero classification.

If $\operatorname{Re}z\ge N+1$, then
$|1-z/j|\ge(N+1-j)/j$ for every $1\le j\le N$, with strict inequality
if $\operatorname{Re}z>N+1$ or $\operatorname{Im}z\ne0$. Multiplication
gives $|G_n(z)|>1$ in the strict cases, proving \eqref{eq:zeros-strip}.

Finally a multiple zero would satisfy both $G_n(\rho)=1$ and
$G_n'(\rho)=0$. By Rolle's theorem all $N-1$ zeros of $G_n'$ are real
and lie in $(1,N)$, where $|G_n|<1$. This is impossible.
\end{proof}

\begin{corollary}[Reflection for odd degree]\label{cor:zeros-reflection}
If $n$ is odd, the multiset of zeros of $P_n^{(1)}$ after removal of
$n+2$ is invariant under $z\mapsto n+2-z$ and $z\mapsto\overline z$.
In particular, the zeros other than $n+2$ occur in reflection orbits
of at most four elements.
\end{corollary}
\begin{proof}
When $N$ is even, \eqref{eq:zeros-reflection} preserves the equation
$G_n(z)=1$. Its solutions $0$ and $N+1$ are exchanged; removing this pair
leaves precisely the stated set. Conjugation preserves the zero set
because all coefficients are real.
\end{proof}

\subsection{Harmonic scaling and the incomplete gamma function}

For a fixed integer $m\ge1$, define the entire function
\begin{equation}\label{eq:zeros-limit-function}
 \mathcal F_m(z)=\frac1{(m-1)!}\int_0^1 t^{m-1}e^{-zt}\,dt
 =\frac{1-e^{-z}\sum_{r=0}^{m-1}z^r/r!}{z^m},
\end{equation}
where the quotient at zero is interpreted as $1/m!$. This function is
 the normalized incomplete gamma function $\gamma^*(m,z)$ in the
 notation of \cite[\S8.2]{DLMF}; see also \cite[\S8.13]{DLMF} and
 \cite{Temme1995} for the classical zero theory of incomplete gamma
 functions. Here the limiting parameter is the interpolation degree,
 not the gamma shape parameter.

\begin{proposition}[Local scaling and zero asymptotics]\label{prop:zeros-scaling}
For fixed $m\ge1$, as $n\to\infty$,
\begin{equation}\label{eq:zeros-scaling}
 \frac{P_n^{(m)}(z/L_n)}{L_n^m}\longrightarrow\mathcal F_m(z)
\end{equation}
locally uniformly in $\mathbb C$. Every zero $w$ of $\mathcal F_m$ is
simple and nonzero. For each such $w$, there is, for all sufficiently
large $n$, a unique zero $\rho_n$ of $P_n^{(m)}$ with $L_n\rho_n$ in a
sufficiently small fixed neighbourhood of $w$. It satisfies
\begin{equation}\label{eq:zeros-general-correction}
 \boxed{\displaystyle \rho_n=\frac{w}{L_n}
 -\frac{\zeta(2)w(w+m-1)}{2L_n^3}+O_{m,w}(L_n^{-4}).}
\end{equation}
Moreover, every finite accumulation point of scaled zeros $L_n\rho_n$
is a zero of $\mathcal F_m$.
\end{proposition}
\begin{proof}
Put $E_j(z)=\sum_{r=0}^j z^r/r!$ for $j\ge0$ and $E_j(z)=0$ for $j<0$.
The logarithmic expansion in \cref{prop:bell} gives, locally uniformly
in $z$,
\[
 G_n(z/L_n)=e^{-z}\left(1-\frac{H_N^{(2)}z^2}{2L_n^2}
 +O(L_n^{-3})\right).
\]
Inverting the expansion and extracting each of the finitely many
coefficients needed in \eqref{eq:beta-series} gives
\[
 \sum_{r=0}^{m-1}\beta_{n,r}(z/L_n)^r
 =E_{m-1}(z)+\frac{H_N^{(2)}z^2}{2L_n^2}E_{m-3}(z)
 +O_m(L_n^{-3}).
\]
These estimates hold on compact subsets of the plane for large $n$,
since $H_N^{(r)}\le\zeta(2)$ for $r\ge2$ and $L_n\to\infty$.
Substitution into \eqref{eq:Pgeneral} shows that
\begin{align}\label{eq:zeros-numerator-expansion}
 z^m\frac{P_n^{(m)}(z/L_n)}{L_n^m}
 ={}&1-e^{-z}E_{m-1}(z)\notag\\
 &+\frac{H_N^{(2)}}{2L_n^2}e^{-z}z^2
     (E_{m-1}(z)-E_{m-3}(z))+O_m(L_n^{-3}).
\end{align}
This proves \eqref{eq:zeros-scaling} away from the origin; uniformity
across the origin follows by the maximum modulus principle on a disc,
since both sides of \eqref{eq:zeros-scaling} are entire.

Write $A_m(z)=1-e^{-z}E_{m-1}(z)$. Then
$A_m'(z)=e^{-z}z^{m-1}/(m-1)!$, and $\mathcal F_m(0)=1/m!$.
Thus all zeros of $\mathcal F_m$ are nonzero and simple. Rouch\'e's
theorem on a sufficiently small disc around $w$ gives a unique scaled
zero $w_n=L_n\rho_n$ tending to $w$. The uniform expansion
\eqref{eq:zeros-numerator-expansion} and $A_m'(w)\ne0$ first imply
$w_n-w=O(L_n^{-2})$ and then
\[
 w_n-w=-\frac{H_N^{(2)}}{2L_n^2}
 \frac{e^{-w}w^2(E_{m-1}(w)-E_{m-3}(w))}{A_m'(w)}
 +O_{m,w}(L_n^{-3}).
\]
The quotient is $w(w+m-1)$, including the case $m=1$. Since
$H_N^{(2)}=\zeta(2)+O(N^{-1})$, division by $L_n$ proves
\eqref{eq:zeros-general-correction}. The last assertion follows directly
from local uniform convergence and continuity.
\end{proof}

In particular the scaled zeros are governed locally by the zeros of a
classical special function. This is a local statement for bounded scaled
variables, not a claim about the global distribution of all $n$ zeros.

\subsection{Higher-order zeta corrections for the first reciprocal power}

For $m=1$, the limit function is $(1-e^{-z})/z$ and its zeros are
$2\pi i\ell$ for $\ell\in\mathbb Z\setminus\{0\}$. The product equation
permits a more precise expansion and reveals the role of higher zeta
values.

\begin{corollary}[Asymptotic zeta corrections]\label{cor:zeros-zeta}
For each fixed $\ell\in\mathbb Z\setminus\{0\}$, let
$\rho_{n,\ell}$ be the zero specified by \cref{prop:zeros-scaling} with
$w=2\pi i\ell$. Then
\begin{equation}\label{eq:zeros-zeta}
 \boxed{\begin{aligned}
 \rho_{n,\ell}={}&\frac{2\pi i\ell}{L_n}
 +\frac{2\pi^2\ell^2\zeta(2)}{L_n^3}
 +\frac{8\pi^3 i\ell^3\zeta(3)}{3L_n^4}\\
 &+O_\ell(L_n^{-5}).
 \end{aligned}}
\end{equation}
In particular,
\[
 \lim_{n\to\infty}L_n^3\operatorname{Re}\rho_{n,\ell}
 =\frac{\pi^4\ell^2}{3},\qquad
 \lim_{n\to\infty}L_n\operatorname{Im}\rho_{n,\ell}=2\pi\ell.
\]
\end{corollary}
\begin{proof}
Near zero the analytic logarithm of $G_n$ with value zero at zero is
$-L_nz-\sum_{r\ge2}H_N^{(r)}z^r/r$. Since $G_n(\rho_{n,\ell})=1$ and
$L_n\rho_{n,\ell}\to w=2\pi i\ell$, the integer in the logarithmic
zero equation is fixed for all sufficiently large $n$, giving
\[
 L_n\rho_{n,\ell}+\sum_{r=2}^\infty
 \frac{H_N^{(r)}}r\rho_{n,\ell}^r=w.
\]
By \eqref{eq:zeros-general-correction}, $\rho_{n,\ell}=w/L_n+O(L_n^{-3})$.
Substitution in this last equation therefore yields
\[
 \rho_{n,\ell}=\frac{w}{L_n}-\frac{H_N^{(2)}w^2}{2L_n^3}
 -\frac{H_N^{(3)}w^3}{3L_n^4}+O_\ell(L_n^{-5}).
\]
For $r=2,3$, the estimate $H_N^{(r)}=\zeta(r)+O(N^{1-r})$
allows replacement by the corresponding zeta values without changing
the error term. Finally insert $w=2\pi i\ell$ and $\zeta(2)=\pi^2/6$.
\end{proof}

More generally, these expansions exist to every fixed order. Indeed,
with $\varepsilon=L_n^{-1}$ consider the analytic implicit equation
\[
 z+\varepsilon\sum_{r=2}^\infty\frac{\zeta(r)}r z^r
 =\varepsilon w.
\]
The analytic implicit function theorem gives a solution near
$(z,\varepsilon)=(0,0)$ with $z(0)=0$. Its coefficients are polynomials
in $w$ and finitely many ordinary zeta values, with rational
coefficients. The difference between the finite harmonic sum in the
exact zero equation and its zeta counterpart is $O(N^{-1}|z|^2)$
uniformly for $|z|\le1/2$, by the integral test. Since the derivative
of the normalized implicit equation is $1+O(\varepsilon^2)$ near the
local solution, this replacement changes the zero by
$O_\ell(N^{-1}L_n^{-3})$, smaller than any fixed inverse power of $L_n$.
Thus the Taylor series of the analytic implicit solution supplies the
full asymptotic expansion of $\rho_{n,\ell}$. This is consistent with
classical analytic inversion (see \cite[\S1.10(vii)]{DLMF}) and the
logarithmic gamma expansion in \cite[\S5.7(i)]{DLMF}.

\subsection{Eventual real-rootedness at fixed degree}

\begin{proposition}[Eventual real-rootedness]\label{prop:zeros-large-m}
Fix $n\ge1$. As $m\to\infty$, all zeros of $P_n^{(m)}$ are eventually
real, positive and simple. They may be labelled $\rho_{n,j}^{(m)}$,
$2\le j\le N$, so that
\begin{equation}\label{eq:zeros-large-m}
 \boxed{\displaystyle
 \rho_{n,j}^{(m)}=j+(-1)^j\frac{(N-1)!}{(j-2)!(N-j)!}
 j^{-m}\bigl(1+O_n(2^{-m})\bigr).}
\end{equation}
A sufficient condition for real-rootedness and simplicity is
\begin{equation}\label{eq:zeros-threshold}
 2^m>(3N-2)(2^{N-1}-1).
\end{equation}
\end{proposition}
\begin{proof}
Let $\ell_j$ be the Lagrange basis polynomial corresponding to node $j$.
Then
\[
 P_n^{(m)}(z)=\sum_{j=1}^N j^{-m}\ell_j(z)
 =\ell_1(z)+O_n(2^{-m}),\qquad
 \ell_1(z)=\prod_{j=2}^N\frac{j-z}{j-1},
\]
with uniform convergence on compact sets, also for all derivatives.
For each $j=2,\ldots,N$, consider the circle $|z-j|=1/3$.
On the union of these circles the Lagrange product formula gives,
for $2\le k\le N$,
\[
 \left|\frac{\ell_k(z)}{\ell_1(z)}\right|
 =\binom{N-1}{k-1}\left|\frac{z-1}{z-k}\right|
 \le(3N-2)\binom{N-1}{k-1}.
\]
Hence $|(P_n^{(m)}-\ell_1)/\ell_1|\le
(3N-2)(2^{N-1}-1)2^{-m}$. Under \eqref{eq:zeros-threshold},
Rouch\'e's theorem gives exactly one zero, counted with multiplicity,
in each of the disjoint discs. Conjugation preserves each disc and
the zero set, so these zeros are real and simple; they are positive
and exhaust the degree $n$.

As $m\to\infty$ these zeros tend to the corresponding $j$. By the real
mean value theorem and $P_n^{(m)}(j)=j^{-m}$, there exists a point
$\xi_m$ between $j$ and $\rho_{n,j}^{(m)}$ such that
\[
 \rho_{n,j}^{(m)}-j=-\frac{j^{-m}}{(P_n^{(m)})'(\xi_m)}.
\]
Since the denominator tends to $\ell_1'(j)\ne0$, the displacement is
$O_n(j^{-m})$. Uniform convergence of derivatives then gives
$(P_n^{(m)})'(\xi_m)=\ell_1'(j)+O_n(2^{-m})$. Finally,
$\ell_1'(j)=(-1)^{j-1}(j-2)!(N-j)!/(N-1)!$ yields
\eqref{eq:zeros-large-m}.
\end{proof}

Thus, if $M_n$ denotes the least positive integer such that all
$P_n^{(m)}$ with $m\ge M_n$ have $n$ simple real zeros, then
\[
 M_n\le1+\left\lfloor\log_2\big((3n+1)(2^n-1)\big)\right\rfloor
 =n+O(\log n).
\]
Here the displayed bound is only sufficient; determining the exact
threshold, the number of real zeros at intermediate exponents, and
simplicity for all $m$ remain questions not settled by the arguments
above.

\section*{Concluding remarks}\label{sec:conclusion}
The reciprocal-power interpolation family admits a common algebraic description through formal inversion of the nodal polynomial. This description links the explicit Bell-polynomial formulas to the recurrences, Newton expansion, and special-function generating functions. The weighted Hurwitz-zeta transform shows how the full interpolation polynomial, rather than just its constant coefficient, enters series representations of integer zeta values. The Dirichlet beta extension preserves the same monomial coefficients and binomial normalization, replacing the zeta weights by an explicit rational convolution. In the zero analysis, half-plane localization holds throughout the family, whereas the harmonic scaling identifies a local limiting zero set and gives corrections expressed through zeta values.

Several questions remain open within the present analysis. These include simplicity of the zeros for arbitrary $m$, the number of real zeros at intermediate values of $m$, and the exact threshold for eventual real-rootedness at fixed degree. A further direction is the global distribution of the zeros as the degree increases, which is not determined by the local scaling limit established here. Further questions include simplifying the beta weights and finding analogous fixed-coefficient transforms for other Dirichlet characters.

\printbibliography

@article{Boyadzhiev2009Euler,
  author       = {Boyadzhiev, Khristo N.},
  title        = {Harmonic Number Identities Via {Euler}'s Transform},
  journal      = {Journal of Integer Sequences},
  year         = {2009},
  volume       = {12},
  eid          = {09.6.1},
  url          = {https://cs.uwaterloo.ca/journals/JIS/VOL12/Boyadzhiev/boyadzhiev3.pdf}
}

@misc{Connon2010Beta,
  author       = {Connon, Donal F.},
  title        = {New rapidly converging series representations for values of the {Riemann} zeta function and the {Dirichlet} beta function},
  year         = {2010},
  eprint       = {1003.4592},
  archivePrefix = {arXiv},
  url          = {https://arxiv.org/abs/1003.4592}
}

@article{HessamiPilehrood2010Beta,
  author       = {Hessami Pilehrood, Khodabakhsh and Hessami Pilehrood, Tatiana},
  title        = {Series acceleration formulas for beta values},
  journal      = {Discrete Mathematics \& Theoretical Computer Science},
  year         = {2010},
  volume       = {12},
  number       = {2},
  pages        = {223--236},
  doi          = {10.46298/dmtcs.504},
  url          = {https://dmtcs.episciences.org/504}
}

@article{Languasco2023,
  author       = {Languasco, Alessandro},
  title        = {A unified strategy to compute some special functions of number-theoretic interest},
  journal      = {Journal of Number Theory},
  year         = {2023},
  volume       = {247},
  pages        = {118--161},
  doi          = {10.1016/j.jnt.2022.12.003},
  eprint       = {2111.07686},
  archivePrefix = {arXiv},
  url          = {https://arxiv.org/abs/2111.07686}
}

@misc{Bhandari2026Beta,
  author       = {Bhandari, Narendra},
  title        = {Positive rational series for reciprocal powers of {Catalan}'s constant and {Dirichlet} beta values},
  year         = {2026},
  eprint       = {2608.00429},
  archivePrefix = {arXiv},
  note         = {Preprint},
  url          = {https://arxiv.org/abs/2608.00429}
}

@article{Sedletskii2004,
  author  = {Sedletskii, A. M.},
  title   = {On the zeros of {Laplace} transforms},
  journal = {Mathematical Notes},
  year    = {2004},
  volume  = {76},
  number  = {6},
  pages   = {824--833},
  doi     = {10.1023/B:MATN.0000049682.65990.e7}
}

@article{Temme1995,
  author  = {Temme, Nico M.},
  title   = {Asymptotics of zeros of incomplete gamma functions},
  journal = {Annals of Numerical Mathematics},
  year    = {1995},
  volume  = {2},
  number  = {1--4},
  pages   = {415--423},
  url     = {https://ir.cwi.nl/pub/2212}
}

@misc{Connon2008,
  author        = {Connon, Donal F.},
  title         = {Euler--{Hurwitz} series and non-linear {Euler} sums},
  year          = {2008},
  eprint        = {0803.1304},
  archivePrefix = {arXiv},
  url           = {https://arxiv.org/abs/0803.1304}
}

@article{Hasse1930,
  author  = {Hasse, Helmut},
  title   = {Ein Summierungsverfahren f\"ur die {Riemannsche} $\zeta$-Reihe},
  journal = {Mathematische Zeitschrift},
  year    = {1930},
  volume  = {32},
  pages   = {458--464},
  doi     = {10.1007/BF01194645}
}

@article{Blagouchine2018,
  author  = {Blagouchine, Iaroslav V.},
  title   = {Three notes on {Ser's} and {Hasse's} representations for the zeta-functions},
  journal = {Integers},
  year    = {2018},
  volume  = {18A},
  eid     = {A3},
  pages   = {1--45},
  eprint  = {1606.02044},
  archivePrefix = {arXiv},
  url     = {https://arxiv.org/abs/1606.02044}
}

@article{FlajoletVepstas2008,
  author  = {Flajolet, Philippe and Vepstas, Linas},
  title   = {On differences of zeta values},
  journal = {Journal of Computational and Applied Mathematics},
  year    = {2008},
  volume  = {220},
  number  = {1--2},
  pages   = {58--73},
  doi     = {10.1016/j.cam.2007.07.040},
  eprint  = {math/0611332},
  archivePrefix = {arXiv},
  url     = {https://arxiv.org/abs/math/0611332}
}

@article{Hoffman1992,
  author  = {Hoffman, Michael E.},
  title   = {Multiple harmonic series},
  journal = {Pacific Journal of Mathematics},
  year    = {1992},
  volume  = {152},
  number  = {2},
  pages   = {275--290},
  url     = {https://msp.org/pjm/1992/152-2/pjm-v152-n2-p04-p.pdf}
}

@article{OhnoOkuda2007,
  author  = {Ohno, Yasuo and Okuda, Jun-Ichi},
  title   = {On the sum formula for the $q$-analogue of non-strict multiple zeta values},
  journal = {Proceedings of the American Mathematical Society},
  year    = {2007},
  volume  = {135},
  number  = {10},
  pages   = {3029--3037},
  doi     = {10.1090/S0002-9939-07-08994-0}
}

@article{XuCai2018,
  author  = {Xu, Ce and Cai, Yulin},
  title   = {On harmonic numbers and nonlinear {Euler} sums},
  journal = {Journal of Mathematical Analysis and Applications},
  year    = {2018},
  volume  = {466},
  number  = {1},
  eprint  = {1609.04924},
  archivePrefix = {arXiv},
  url     = {https://arxiv.org/abs/1609.04924}
}

@misc{DLMF,
  author = {{National Institute of Standards and Technology}},
  title  = {{NIST} Digital Library of Mathematical Functions},
  note   = {Chapters 5 and 25},
  url    = {https://dlmf.nist.gov/}
}

@article{Ser1926,
  author  = {Ser, Joseph},
  title   = {Sur une expression de la fonction $\zeta(s)$ de {Riemann}},
  journal = {Comptes rendus de l'Académie des sciences de Paris},
  year    = {1926},
  volume  = {182},
  pages   = {1075--1077}
}

@article{Sondow1994,
  author  = {Sondow, Jonathan},
  title   = {Analytic continuation of {Riemann's} zeta function and values at negative integers via {Euler's} transformation of series},
  journal = {Proceedings of the American Mathematical Society},
  year    = {1994},
  volume  = {120},
  number  = {2},
  pages   = {421--424}
}

@book{Comtet1974,
  author    = {Comtet, Louis},
  title     = {Advanced Combinatorics: The Art of Finite and Infinite Expansions},
  publisher = {D. Reidel Publishing Company},
  address   = {Dordrecht},
  year      = {1974},
  doi       = {10.1007/978-94-010-2196-8}
}

@article{Leung2018,
  author  = {Leung, Ho-Hon},
  title   = {Another Identity for Complete Bell Polynomials Based on {Ramanujan}'s Congruences},
  journal = {Journal of Integer Sequences},
  volume  = {21},
  year    = {2018},
  eid     = {18.6.4},
  note    = {Article 18.6.4},
  url     = {https://cs.uwaterloo.ca/journals/JIS/VOL21/Leung/leung4.pdf}
}

@book{Macdonald1995,
  author    = {Macdonald, I. G.},
  title     = {Symmetric Functions and Hall Polynomials},
  edition   = {2},
  publisher = {Oxford University Press},
  address   = {Oxford},
  year      = {1995},
  note      = {See Chapter I, Section 2},
  url       = {https://academic.oup.com/book/52932}
}

@book{Stanley1999,
  author    = {Stanley, Richard P.},
  title     = {Enumerative Combinatorics},
  volume    = {2},
  publisher = {Cambridge University Press},
  address   = {Cambridge},
  year      = {1999},
  note      = {See Chapter 5 (the exponential formula)},
  url       = {https://www.cambridge.org/core/books/enumerative-combinatorics/D8DDDFF7E8EBF0BCFE99F5E6918CE2A8}
}
\end{document}